\documentclass{amsart}

\usepackage{amsmath}
\usepackage{amssymb}

\newtheorem{theorem}[subsection]{Theorem}
\newtheorem{proposition}[subsection]{Proposition}
\newtheorem{corollary}[subsection]{Corollary}
\newtheorem{lemma}[subsection]{Lemma}

\newtheorem{claim}{Claim}

\newtheorem{question}[subsection]{Question}

\providecommand{\supp}{\mathop{\rm supp}\nolimits}
\providecommand{\cod}{\mathop{\rm cod}\nolimits}

\newcommand{\wh}{\widehat}

\begin{document}

\title[Additive structures in sumsets]
  {Additive structures in sumsets}

\author{TOM SANDERS}
\address{Department of Pure Mathematics and Mathematical Statistics\\
University of Cambridge\\
Wilberforce Road\\
Cambridge CB3 0WA\\
England } \email{t.sanders@dpmms.cam.ac.uk}

\begin{abstract}
Suppose that $A$ and $A'$ are subsets of $\mathbb{Z}/N\mathbb{Z}$.
We write $A+A'$ for the set $\{a+a':a \in A \textrm{ and } a' \in
A'\}$ and call it the \emph{sumset} of $A$ and $A'$. In this paper
we address the following question. Suppose that $A_1,...,A_m$ are
subsets of $\mathbb{Z}/N\mathbb{Z}$. Does $A_1+...+A_m$ contain a
long arithmetic progression?

The situation for $m=2$ is rather different from that for $m \geq
3$. In the former case we provide a new proof of a result due to
Green. He proved that $A_1+A_2$ contains an arithmetic progression
of length roughly $\exp (c\sqrt{\alpha_1\alpha_2 \log N})$ where
$\alpha_1$ and $\alpha_2$ are the respective densities of $A_1$
and $A_2$. In the latter case we improve the existing estimates.
For example we show that if $A \subset \mathbb{Z}/N\mathbb{Z}$ has
density $\alpha \gg \sqrt{\log \log N/\log N}$ then $A+A+A$
contains an arithmetic progression of length $N^{c\alpha}$. This
compares with the previous best of $N^{c\alpha^{2+\varepsilon}}$.

Two main ingredients have gone into the paper. The first is the
observation that one can apply the iterative method to these
problems using some machinery of Bourgain. The second is that we
can localize a result due to Chang regarding the large spectrum of
$L^2$-functions. This localization seems to be of interest in its
own right and has already found one application elsewhere.
\end{abstract}

\maketitle

\section{Introduction}

As indicated in the abstract we are interested in the following
question.
\begin{question}
Suppose that $A_1,...,A_m \subset \mathbb{Z}/N\mathbb{Z}$. Does
$A_1+...+A_m$ contain a long arithmetic progression?
\end{question}
The case $m=2$ is much harder than $m \geq 3$; the best known
bounds lie with Green, \cite{BJGAA}, who proved the following
result.
\begin{theorem}\label{green1}
Suppose that $A_1,A_2 \subset \mathbb{Z}/N\mathbb{Z}$. Suppose
that $\alpha$, the geometric mean of the densities of the $A_i$s,
is positive. Then $A_1+A_2$ contains an arithmetic progression of
length at least $\exp (c((\alpha^2\log N)^{\frac{1}{2}} - \log
\log N))$ for some absolute constant $c>0$.
\end{theorem}
The next result is proved for $m=3$ and $A_1=A_2=A_3$ in
\cite{BJGAA}, although the more general conclusion can easily be
read out of that paper. The special case $m=4$, $A_1=A_2=A$,
$A_3=A_4=-A$ can be read out of earlier work of Chang \cite{MCC}.
\begin{theorem}\label{green2}
Suppose that $m \geq 3$ and $A_1,...,A_m \subset
\mathbb{Z}/N\mathbb{Z}$. Suppose that $\alpha$, the geometric mean
of the densities of the $A_i$s, is positive. Then $A_1+...+A_m$
contains an arithmetic progression of length at least $c\alpha^C
N^{cm^{-1}\alpha^{\frac{2}{m-2}}(\log \alpha^{-1})^{-1}}$ for some
absolute constants $C,c>0$.
\end{theorem}
In this paper we prove the following two results.
\begin{theorem}\label{green1new}
Suppose that $A_1,A_2 \subset \mathbb{Z}/N\mathbb{Z}$. Suppose
that $\alpha$, the geometric mean of the densities of the $A_i$s,
is positive. Then $A_1+A_2$ contains an arithmetic progression of
length at least $\exp c(((\alpha^2\log N)^{\frac{1}{2}} - \log
\alpha^{-1}\log \log N)$ for some absolute constant $c>0$.
\end{theorem}
\begin{theorem}\label{greennew}
Suppose that $m \geq 3$ and $A_1,...,A_m \subset
\mathbb{Z}/N\mathbb{Z}$. Write $\alpha$ for the geometric mean of
the densities of the $A_i$s. Then $A_1+...+A_m$ contains an
arithmetic progression of length
$c\alpha^{Cm^3\alpha^{-\frac{1}{m-2}}}
N^{cm^{-2}\alpha^{\frac{1}{m-2}}}$ for some absolute constants
$C,c>0$.
\end{theorem}
Theorem \ref{greennew} is stronger than Theorem \ref{green2}, while
Theorem \ref{green1new} is marginally weaker than Theorem
\ref{green1}. Despite this, we believe that merit can still be found
in our approach to Theorem \ref{green1new} for two reasons:
\begin{itemize}
\item The real strength of Theorems \ref{green1} and
\ref{green1new} is when $A_1$ and $A_2$ are thick sets, and in that
case the slightly weaker error term in Theorem \ref{green1new} plays
no part. The case when $A_1$ and $A_2$ are thin is addressed by
Croot, Ruzsa and Schoen in \cite{ESCIZRTS}. \item Green's proof of Theorem \ref{green1} is a
tour de force combining a number of powerful analytic tools in a
highly non-trivial way. By contrast our method is conceptually
simpler although probably technically more challenging.
\end{itemize}
In any case the merit of our approach can, perhaps, be most easily
seen in the finite field setting where the two methods give the
same result and the technicalities in our approach disappear to
leave a fairly simple argument.

To put our refinement of Theorem \ref{green2} in context take, for
example, $m=3$ and $A_1=A_2=A_3=A$. Theorem \ref{greennew} is then
equivalent to the fact that there are absolute constants $C,c>0$
such that if
\begin{equation*}
\alpha \geq C\sqrt{\frac{\log \log N}{\log N}} \textrm{ then }
A+A+A \textrm{ contains a progression of length } N^{c\alpha}.
\end{equation*}
The previous best is essentially equivalent to the existence of
absolute constants $C,c>0$ such that if
\begin{equation*}
\alpha \geq C\sqrt{\frac{(\log \log N)^2}{\log N}} \textrm{ then }
A+A+A \textrm{ contains a progression of length }
N^{c\alpha^{2+\varepsilon}}.
\end{equation*}
Here, of course, $\alpha^{2+\varepsilon}$ is shorthand for
$\alpha^2$ up to some logarithmic factors.

It is also the case, as we shall see in the next section, that the
proof behind Theorem \ref{greennew} gives stronger structural
information than Theorem \ref{green2}.

\section{The Fourier transform, Bohr neighborhoods and additive structure}

In this section we develop the results behind Theorems
\ref{green1new} and \ref{greennew}, and identify some of the
mathematics necessary to prove them. Our main tool is the Fourier
transform; we take a moment to set our notation in this regard.

Suppose that $G$ is a compact Abelian group. We write $\wh{G}$ for
the dual group, that is the discrete Abelian group of continuous
homomorphisms $\gamma:G \rightarrow S^1$, where $S^1:=\{z \in
\mathbb{C}:|z|=1\}$. Although the natural group operation on
$\wh{G}$ corresponds to pointwise multiplication of characters we
shall denote it by `$+$' in alignment with contemporary work. $G$
may be endowed with Haar measure $\mu_G$ normalized so that
$\mu_G(G)=1$ and as a consequence we may define the Fourier
transform $\widehat{.}:L^1(G) \rightarrow \ell^\infty(\widehat{G})$
which takes $f \in L^1(G)$ to
\begin{equation*}
\widehat{f}: \widehat{G} \rightarrow \mathbb{C}; \gamma \mapsto
\int_{x \in G}{f(x)\overline{\gamma(x)}d\mu_G(x)}.
\end{equation*}

We can define a natural valuation on $S^1$, namely, given $z \in
S^1$ we let $\|z\|:=|\theta|$ where $\theta$ is the unique element
of the interval $(-1/2,1/2]$ such that $x=\exp(2\pi i\theta)$. This
valuation can be used to measure how far $\gamma(x)$ is from 1.
Suppose that $\Gamma \subset \wh{G}$ and $\delta \in (0,1]$. We
define
\begin{equation*}
B(\Gamma,\delta):=\{x \in G: \|\gamma(x)\| \leq \delta \textrm{ for
all }\gamma \in \Gamma\},
\end{equation*}
and call such a set a \emph{Bohr set} and a translate of such a set
a \emph{Bohr neighborhood}. The following is an easy pigeonhole
argument and gives an estimate for the volume of these sets. See
Lemma 4.20 in \cite{TCTVHV} for the details.
\begin{lemma}\label{bohrsize}
Suppose $G$ is a compact Abelian group, $\Gamma$ a set of $d$
characters on $G$ and $\delta \in (0,1]$. Then
$\mu_G(B(\Gamma,\delta)) \geq \delta^d$.
\end{lemma}
Consequently we can write $\beta_{\Gamma,\delta}$, or simply $\beta$
or $\beta_\delta$ if the parameters are implicit, for the measure
induced on $B(\Gamma,\delta)$ by $\mu_G$ and normalized so that
$\|\beta_{\Gamma,\delta}\|_1=1$. This is sometimes referred to as
the \emph{normalized Bohr cutoff}. We write $\beta'$ for
$\beta_{\Gamma',\delta'}$, or $\beta_{\Gamma,\delta'}$ if no
$\Gamma'$ has been defined. We have a similar convention for
$\beta''$.

In \S\ref{lfanal} we shall see that Bohr sets in fact behave as sort
of approximate groups and in particular they are highly additively
structured. When $G$ is a cyclic group this translates to Bohr sets
containing long arithmetic progressions; specifically the following
is another easy application of the pigeonhole principle. Again, see
\cite{TCTVHV} for details.
\begin{lemma}\label{triviallemma}
Suppose that $G=\mathbb{Z}/N\mathbb{Z}$. Suppose that
$B(\Gamma,\delta)$ is a Bohr set with $|\Gamma|=d$. Then
$B(\Gamma,\delta)$ contains an arithmetic progression of length
$\delta N^{\frac{1}{d}}$.
\end{lemma}

Theorem \ref{green2} follows immediately from Lemma
\ref{triviallemma} and the following result about Bohr
neighborhoods.
\begin{theorem}\label{green2'}
Suppose that $G$ is a compact Abelian group. Suppose that $m \geq
3$ and $A_1,...,A_m \subset G$. Suppose that $\alpha$, the
geometric mean of the densities of the $A_i$s, is positive. Then
$A_1+...+A_m$ contains a translate of a Bohr set
$B(\Gamma,\delta)$ with
\begin{equation*}
|\Gamma| \ll m\alpha^{-\frac{2}{m-2}}\log \alpha^{-1} \textrm{ and
} \log \delta^{-1} \ll \log \alpha^{-1}.
\end{equation*}
\end{theorem}
In this paper we prove the following refinement.
\begin{theorem}\label{aaageneral}
Suppose that $G$ is a compact Abelian group. Suppose that $m \geq
3$ and $A_1,...,A_m \subset G$. Suppose that $\alpha$, the
geometric mean of the densities of the $A_i$s, is positive. Then
$A_1+...+A_m$ contains an Bohr neighborhood $B(\Gamma,\delta)$
with
\begin{equation*}
|\Gamma| \ll m^2\alpha^{-\frac{1}{m-2}} \textrm{ and } \log
\delta^{-1} \ll m^3\alpha^{-\frac{1}{m-2}}\log \alpha^{-1}.
\end{equation*}
\end{theorem}
The dimension of our Bohr set is much smaller than that found by
Green and it is this which ensures that the arithmetic progression
we find (when there is one at all) is much longer.

Although we do not require it, there is an important strengthening
of Lemma \ref{triviallemma} due to Ruzsa, \cite {IZRF}. If
$M=P_1+...+P_d$ where $P_1,...,P_d$ are arithmetic progressions
then we call $M$ a \emph{$d$-dimensional generalized arithmetic
progression}.
\begin{lemma}
Suppose that $G=\mathbb{Z}/N\mathbb{Z}$. Suppose that
$B(\Gamma,\delta)$ is a Bohr set with $|\Gamma|=d$. Then
$B(\Gamma,\delta)$ contains a $d$-dimensional generalized
arithmetic progression of size at least $(\delta/d)^d N$.
\end{lemma}
Typically this generalized progression occupies a large (roughly
$(Cd)^{-d}$) proportion of the Bohr set, and it can be instructive
to think of Bohr neighborhoods as generalized progressions. In our
results, then, we could use this lemma to draw the stronger
conclusion that $m$-fold sumsets (for $m \geq 3$) contain large
multidimensional progressions, but we believe that the results are
most easily digested in the form stated.

The paper now splits into five further sections. In the next three
we develop the necessary tools for analyzing functions on Bohr sets.
The first of these presents the basics, the second establishes our
new version of Chang's theorem relative to Bohr sets, and the third
recalls some standard density increment lemmas. The last two
sections of the paper prove the results we have promised.

\section{Local Fourier analysis on compact Abelian
groups}\label{lfanal}

Given $f \in L^1(G)$ we often want to approximate $f$ by a less
complicated function. One way to do this is to approximate $f$ by
its expectation on approximate level sets of characters i.e. sets on
which characters do not vary too much. To analyze the error in doing
this we restrict the function to these approximate level sets and
use the Fourier transform on the restricted function. Specifically,
if $\Gamma$ is a set of characters and $x' + \Gamma^\perp$ (a
maximal joint level set of the characters in $\Gamma$) has positive
measure in $G$ then it is easy to localize the Fourier transform to
$x'+\Gamma^{\perp}$:
\begin{equation*}
L^1(x'+\mu_{\Gamma^\perp}) \rightarrow \ell^\infty(\wh{G}); f
\mapsto \wh{fd(x'+\mu_{\Gamma^\perp})}.
\end{equation*}
Note that the right hand side is constant on cosets of
$\Gamma^{\perp\perp}$ and so is really an element of
$\ell^\infty(\wh{G}/\Gamma^{\perp\perp})$.

Bourgain, in \cite{JB}, observed that one can localize the Fourier
transform to typical approximate level sets and retain approximate
versions of a number of the standard results for the Fourier
transform on compact Abelian groups. Since his original work
various expositions and extensions of the work have appeared most
notably in the various papers of Green and Tao. Indeed all the
results of this section can be found in \cite{BJGTCTU3}, for example.

Annihilators are subgroups of $G$; a Bohr set is a sort of
approximate annihilator and, consequently, we would like it to
behave like a sort of approximate subgroup. Suppose that $\eta \in
(0,1]$. Then $B({\Gamma},\delta)+B({\Gamma},\eta\delta) \subset
B({\Gamma},(1+\eta)\delta)$. If $B({\Gamma},(1+\eta)\delta)$ is not
much bigger than $B({\Gamma},\delta)$ then we have a sort of
approximate additive closure in the sense that
$B({\Gamma},\delta)+B({\Gamma},\eta\delta) \approx
B({\Gamma},(1+\eta)\delta)$. Not all Bohr sets have this property.
However, Bourgain showed that typically they do. For our purposes we
have the following proposition.
\begin{proposition}\label{ubreg}
Suppose that $G$ is a compact Abelian group, $\Gamma$ a set of $d$
characters on $G$ and $\delta \in (0,1]$. There is an absolute
constant $c_{\mathcal{R}}>0$ and a $\delta' \in [\delta/2,\delta)$
such that
\begin{equation}\label{reg}
\frac{\mu_G(B(\Gamma,(1+\kappa)\delta'))}{\mu_G(B(\Gamma,\delta'))}
= 1 + O(|\kappa|d)
\end{equation}
whenever $|\kappa|d \leq c_{\mathcal{R}}$.
\end{proposition}
This result is not as easy as the rest of the section, it uses a
covering argument; a nice proof can be found in \cite{BJGTCTU3}. We
say that $\delta'$ is \emph{regular for $\Gamma$} or that
$B(\Gamma,\delta')$ is \emph{regular} if
\begin{equation*}
\frac{\mu_G(B(\Gamma,(1+\kappa)\delta'))}{\mu_G(B(\Gamma,\delta'))}
= 1 + O(|\kappa|d) \textrm{ whenever } |\kappa|d \leq
c_{\mathcal{R}}.
\end{equation*}
It is regular Bohr sets to which we localize the Fourier transform.
We require a little more notation regarding measures. As usual if
$X$ is a topological space we write $M(X)$ for the regular
complex-valued Borel measures on $X$. If $\mu \in M(G)$ then $\supp
\mu$ denotes the support of $\mu$, and if $x \in G$ as well then
$x+\mu$ denotes the measure $\mu$ translated by $x$.

We begin by observing that normalized regular Bohr cutoffs are
approximately translation invariant and so function as normalized
approximate Haar measures.
\begin{lemma}
\label{contlem}\emph{(Normalized approximate Haar measure)} Suppose
that $G$ is a compact Abelian group and $B(\Gamma,\delta)$ is a
regular Bohr set. If $y \in B(\Gamma,\delta')$ then
$\|(y+\beta_\delta) - \beta_\delta\| \ll d\delta'\delta^{-1}$.
\end{lemma}
The proof follows immediately from the definition of regularity.
In applications the following two simple corollaries will be
useful but they should be ignored until they are used.
\begin{corollary}
\label{contlemcor} Suppose that $G$ is a compact Abelian group and
$B(\Gamma,\delta)$ is a regular Bohr set. If $\mu \in
M(B(\Gamma,\delta'))$ then $\|\beta\ast \mu - \beta\int{d\mu}\|
\ll \|\mu\|d\delta'\delta^{-1}$.
\end{corollary}
\begin{corollary}
\label{contlemcor2} Suppose that $G$ is a compact Abelian group
and $B(\Gamma,\delta)$ is a regular Bohr set. Suppose that $f \in
L^\infty(G)$. If $x-y \in B(\Gamma,\delta')$ then $|f \ast
\beta(x) - f \ast \beta(y)| \ll \|f\|_\infty d\delta'\delta^{-1}$.
\end{corollary}

With an approximate Haar measure we are in a position to define
the local Fourier transform: Suppose that $\Gamma$ is a finite set
of characters, $\delta$ is regular for $\Gamma$ and $x' \in G$.
Then we define the Fourier transform local to
$x'+B(\Gamma,\delta)$ by
\begin{equation*}
L^1(x'+\beta_{\Gamma,\delta}) \rightarrow \ell^\infty(\wh{G}); f
\mapsto \wh{fd(x'+\beta_{\Gamma,\delta})},
\end{equation*}
where we take the convention that
\begin{equation*}
L^1(\mu):=\{f \in L^1(G):\supp f \subset \supp \mu \textrm{ and }
\int{|f|d\mu} < \infty \}.
\end{equation*}
The translation of the Bohr set by $x'$ simply twists the Fourier
transform and is unimportant for the most part so we tend to
restrict ourselves to the case when $x'=0$.

$\wh{fd\mu_{\Gamma^\perp}}$ was constant on cosets of
$\Gamma^{\perp\perp}$. In the approximate setting have an
approximate analogue of this. First the analogue of
$\Gamma^{\perp\perp}$; there are a number of possibilities:
\begin{eqnarray*}
&\{\gamma:|1-\gamma(x)| \leq \epsilon \textrm{ for all
} x \in B(\Gamma,\delta)\} &\textrm{ for } \epsilon \in (0,1] \\
&\{\gamma:|1-\wh{\beta}(\gamma)| \leq \epsilon\}
&\textrm{ for } \epsilon  \in (0,1]\\
&\{\gamma:|\wh{\beta}(\gamma)| \geq \epsilon \} &\textrm{ for }
\epsilon  \in (0,1].
\end{eqnarray*}
In applications each of these classes of sets is useful and so we
should like all of them to be approximately equivalent. There is a
clear chain of inclusions between the classes:
\begin{equation*}
\{\gamma:|1-\gamma(x)| \leq \epsilon \textrm{ for all } x \in
B(\Gamma,\delta)\} \subset \{\gamma:|1-\wh{\beta}(\gamma)| \leq
\epsilon\} \subset \{\gamma:|\wh{\beta}(\gamma)| \geq 1-\epsilon
\}
\end{equation*}
for $\epsilon \in (0,1]$. For a small cost in the width of the
Bohr set we can ensure that the sets in the third class are
contained in those in the first.
\begin{lemma} \label{nestsupport}
Suppose that $G$ is a compact Abelian group and $B(\Gamma,\delta)$
is a regular Bohr set. Suppose that $\eta_1,\eta_2>0$. Then there
is a $\delta' \gg \eta_1\eta_2\delta/d$ such that
\begin{equation*}
\{\gamma: |\wh{\beta}(\gamma)| \geq \eta_1\} \subset \{\gamma:
|1-\gamma(x)| \leq \eta_2 \textrm{ for all } x \in
B(\Gamma,\delta')\}.
\end{equation*}
\end{lemma}
The lemma follows easily from Lemma \ref{contlem}.

\section{The structure of sets of characters supporting large
values of the local Fourier transform}

Green and Tao in \cite{BJGTCTU3} were the first to prove the
following proposition which captures a version of Bessel's
inequality local to Bohr sets in a form useful for applications.
\begin{proposition}\label{local Bessels bound} Suppose that $G$ is a compact
Abelian group and $B(\Gamma,\delta)$ is a regular Bohr set.
Suppose that $f \in L^2(\beta)$ and $\epsilon,\eta \in (0,1]$.
Then there is a set of characters $\Lambda$ and a $\delta' \in
(0,1]$ such that
\begin{equation*}
|\Lambda| \ll
\epsilon^{-2}\|f\|_{L^1(\beta)}^{-2}\|f\|_{L^2(\beta)}^2 \textrm{
and } \delta' \gg \eta\delta/d,
\end{equation*}
and
\begin{equation*}
\{\gamma \in \wh{G}:|\wh{fd\beta}(\gamma)| \geq \epsilon
\|f\|_{L^1(\beta)}\} \subset \{\gamma \in \wh{G}:|1-\gamma(x)|
\leq \eta \textrm{ for all } x \in B(\Gamma \cup
\Lambda,\delta')\}.
\end{equation*}
\end{proposition}
It is instructive to consider this result in the case when
$B(\Gamma,\delta)=G$. We are then given a set of characters
$\Lambda$ and a $\delta' \in (0,1]$ such that
\begin{equation*}
|\Lambda| \ll
\epsilon^{-2}\|f\|_{L^1(\mu_G)}^{-2}\|f\|_{L^2(\mu_G)}^2 \textrm{
and } \delta' \gg \eta,
\end{equation*}
and
\begin{equation*}
  \{\gamma \in \wh{G}:|\wh{f}(\gamma)| \geq \epsilon
\|f\|_{L^1(\mu_G)}\} \subset \{\gamma \in \wh{G}:|1-\gamma(x)| \leq
\eta \textrm{ for all } x \in B(\Lambda,\delta')\}.
\end{equation*}
Now it is natural to ask what this has to do with Bessel's
inequality. Write $\Gamma:=\{\gamma \in \wh{G}:|\wh{f}(\gamma)| \geq
\epsilon \|f\|_{L^1(\mu_G)}\}$. Then
\begin{equation*}
  |\Gamma|.\epsilon^2\|f\|_{L^1(\mu_G)}^2 \leq \sum_{\gamma \in
  \Gamma}{|\wh{f}(\gamma)|^2} \leq \sum_{\gamma \in
  \wh{G}}{|\wh{f}(\gamma)|^2} \leq \|f\|_{L^2(\mu_G)}^2
\end{equation*}
by Bessel's inequality. It follows that
\begin{equation*}
  |\Gamma| \leq
  \epsilon^{-2}\|f\|_{L^1(\mu_G)}^{-2}\|f\|_{L^2(\mu_G)}^2,
\end{equation*}
and, moreover, it is easy to check that there is a $\delta' \gg
\eta$ such that
\begin{equation*}
  \Gamma \subset \{\gamma \in \wh{G}:|1-\gamma(x)| \leq \eta
  \textrm{ for all } x \in B(\Gamma,\delta')\}.
\end{equation*}
Setting $\Lambda:=\Gamma$ yields the conclusion of Proposition
\ref{local Bessels bound} in this special case. Restricting
functions to Bohr sets complicates matters. However, there are some
easy rules of thumb to bear in mind. The bound on $|\Lambda|$ is
very important. The dependence of $\delta'$ on $\delta$ needs to be
linear because we intend to iterate the procedure, however, the
ratio $\delta^{-1}\delta'$ can be very much smaller without any
tangible cost.

In the next proposition we trade a worse ratio $\delta^{-1}\delta'$,
which has little impact, for a large improvement in the bound on
$|\Lambda|$.
\begin{proposition}\label{local Changs bound} Suppose that $G$ is a compact
Abelian group and $B(\Gamma,\delta)$ is a regular Bohr set.
Suppose that $f \in L^2(\beta)$ and $\epsilon,\eta \in (0,1]$.
Then there is a set of characters $\Lambda$ and a $\delta' \in
(0,1]$ such that
\begin{equation*}
|\Lambda| \ll \epsilon^{-2}\log
\|f\|_{L^1(\beta)}^{-2}\|f\|_{L^2(\beta)}^2 \textrm{ and } \delta'
\gg \delta\eta\epsilon^2/d^2 \log
\|f\|_{L^1(\beta)}^{-2}\|f\|_{L^2(\beta)}^2,
\end{equation*}
and
\begin{equation*}
\{\gamma \in \wh{G}:|\wh{fd\beta}(\gamma)| \geq \epsilon
\|f\|_{L^1(\beta)}\} \subset \{\gamma \in \wh{G}:|1-\gamma(x)|
\leq \eta \textrm{ for all } x \in B(\Gamma \cup
\Lambda,\delta')\}.
\end{equation*}
\end{proposition}
This result can be seen as a local version of Chang's theorem
(from \cite{MCC}); indeed, the proof is essentially a combination
of the ideas used to prove that theorem and those used to prove
Proposition \ref{local Bessels bound}. The key tool in Chang's
theorem is that of dissociativity; in the local setting we use the
following version of it. If $\Lambda$ is a set of characters on
$G$ and $m:\Lambda \rightarrow \mathbb{Z}$ has finite support then
put
\begin{equation*}
m.\Lambda:=\sum_{\lambda \in \Lambda}{m_\lambda.\lambda} \textrm{
and } |m|:=\sum_{\lambda \in \Lambda}{|m_\lambda|},
\end{equation*}
where the second `$.$' in the first definition denotes the natural
action of $\mathbb{Z}$ on $\wh{G}$. If $S$ is a non-empty symmetric
neighborhood of $0_{\wh{G}}$ then we say that $\Lambda$ is
\emph{$S$-dissociated} if
\begin{equation*}
m.\Lambda \in S \Rightarrow m\equiv 0.
\end{equation*}
The usual definition of dissociativity corresponds to taking
$S=\{0_{\wh{G}}\}$.

Proposition \ref{local Changs bound} follows straightforwardly from
the next two lemmas.
\begin{lemma}\label{spanning chang}
Suppose that $G$ is a compact Abelian group and $B(\Gamma,\delta)$
is a regular Bohr set. Suppose that $\eta',\eta \in (0,1]$ and
$\Delta$ is a set of characters on $G$. If $\Lambda$ is a maximal
$\{\gamma:|\wh{\beta}(\gamma)| \geq \eta' \}$-dissociated subset
of $\Delta$ then there is a $\delta' \gg
\min\{\eta/|\Lambda|,\eta' \eta \delta/d\}$ such that
\begin{equation*}
\Delta \subset \{\gamma:|1-\gamma(x)| \leq \eta \textrm{ for all }
x \in B(\Gamma \cup \Lambda,\delta')\}.
\end{equation*}
\end{lemma}
\begin{lemma}\label{content chang}
Suppose that $G$ is a compact Abelian group and $B(\Gamma,\delta)$
is a regular Bohr set. Suppose that $0 \not \equiv f \in
L^2(\beta)$, $k \in \mathbb{N}$ and $\epsilon,\eta \in (0,1]$. Then
there is a $\delta' \gg \delta/dk$ regular for $\Gamma$ such that if
$\Lambda$ is a $\{\gamma:|\wh{\beta'}(\gamma)| \geq
1/3\}$-dissociated subset of $\{\gamma \in
\wh{G}:|\wh{fd\beta}(\gamma)| \geq \epsilon \|f\|_{L^1(\beta)}\}$
with size at most $k$, then
\begin{equation*}
|\Lambda| \ll \epsilon^{-2} \log
\|f\|_{L^1(\beta)}^{-2}\|f\|_{L^2(\beta)}^2.
\end{equation*}
\end{lemma}

\subsection{Proof of Lemma \ref{spanning chang}}

The lemma is really rather simple to prove although technical. It
rests on localizing the following simple observation of Chang
\cite{MCC}.
\begin{lemma}
Suppose that $G$ is a compact Abelian group. Suppose that $\Delta$
is a set of characters on $G$ and $\Lambda$ is a maximal
dissociated subset of $\Delta$. Then $\Delta \subset \langle
\Lambda \rangle$.
\end{lemma}
Here $\langle \Lambda \rangle$ denotes all the finite $\pm$-sums
of elements of $\Lambda$ i.e.
\begin{equation*}
\langle \Lambda \rangle:=\{m.\Lambda: m :\Lambda \rightarrow
\{-1,0,1\} \textrm{ and } |\supp m| < \infty\}.
\end{equation*}
For our purposes we have the following.
\begin{lemma}
Suppose that $G$ is a compact Abelian group. Suppose that $S$ is a
non-empty symmetric neighborhood of $0_{\wh{G}}$. Suppose that
$\Delta$ is a set of characters on $G$ and $\Lambda$ is a maximal
$S$-dissociated subset of $\Delta$. Then $\Delta \subset \langle
\Lambda \rangle + S$.
\end{lemma}
\begin{proof}
If $\lambda_0 \in \Delta \setminus (\langle \Lambda \rangle +S)$
then we put $\Lambda':=\Lambda \cup \{ \lambda_0 \}$, which is a
strict superset of $\Lambda$, and a subset of $\Delta$. It turns
out that $\Lambda'$ is also $S$-dissociated which contradicts the
maximality of $\Lambda$. Suppose that $m:\Lambda' \rightarrow
\{-1,0,1\}$ and $m.\Lambda' \in S$. Then we have three
possibilities for the value of $m_{\lambda_0}$:
\begin{enumerate}
\item $m.\Lambda' = \lambda_0 + m|_\Lambda.\Lambda$, in which case
$\lambda_0 \in -m|_\Lambda.\Lambda + S \subset \langle \Lambda
\rangle + S$ - a contradiction; \item $m.\Lambda' = - \lambda_0 +
m|_\Lambda.\Lambda$, in which case $\lambda_0 \in
m|_\Lambda.\Lambda - S \subset \langle \Lambda \rangle + S$ - a
contradiction; \item $m.\Lambda' = m|_\Lambda.\Lambda$, in which
case $m|_\Lambda\equiv 0$ since $\Lambda$ is $S$-dissociated and
hence $m \equiv 0$.
\end{enumerate}
It follows that $m.\Lambda' \in S\Rightarrow m\equiv 0$ i.e.
$\Lambda'$ is $S$-dissociated as claimed.
\end{proof}
Lemma \ref{spanning chang} then follows from the above and the
following lemma.
\begin{lemma}
Suppose that $G$ is a compact Abelian group and $B(\Gamma,\delta)$
is a regular Bohr set. Suppose that $\eta',\eta \in (0,1]$ and
$\Lambda$ is a set of characters on $G$. Then there is a $\delta'
\gg \min\{\eta/|\Lambda|,\eta'\eta\delta/d\}$ such that
\begin{equation*}
\langle \Lambda \rangle + \{\gamma:|\wh{\beta}(\gamma)| \geq \eta'
\} \subset \{\gamma:|1-\gamma(x)| \leq \eta \textrm{ for all } x
\in B(\Gamma \cup \Lambda,\delta')\}.
\end{equation*}
\end{lemma}
\begin{proof}
The lemma has two parts.
\begin{enumerate}
\item If $\lambda \in \langle \Lambda \rangle$ then
\begin{equation*}
|1-\lambda(x)| \leq \sum_{\lambda' \in \Lambda}{|1-\lambda'(x)|},
\end{equation*}
so there is a $\delta'' \gg \eta/|\Lambda|$ such that
\begin{equation*}
\langle \Lambda \rangle \subset \{\gamma: |1-\gamma(x)| \leq
\eta/2 \textrm{ for all } x \in B(\Lambda,\delta'')\}.
\end{equation*}
\item By Lemma \ref{nestsupport} there is a $\delta''' \gg
\eta\eta'\delta/d$ such that
\begin{equation*}
\{\gamma:|\wh{\beta}(\gamma)| \geq \eta' \} \subset
\{\gamma:|1-\gamma(x)| \leq \eta/2 \textrm{ for all } x \in
B(\Gamma,\delta''')\}.
\end{equation*}
\end{enumerate}
Taking $\delta' = \min\{\delta'',\delta'''\}$ we have the result
by the triangle inequality.
\end{proof}

\subsection{Proof of Lemma \ref{content chang}}

The proof has three main ingredients.

\begin{itemize}
\item \emph{(Rudin's inequality)} In \cite{MCC} Chang uses the
following dual statement of Rudin's inequality to prove her
theorem.
\begin{proposition}\label{dualrudin}
Suppose that $G$ is a compact Abelian group. If $\Lambda$ is
dissociated then
\begin{equation*}
\|\wh{f}|_\Lambda\|_2 \ll \sqrt{p}\|f\|_{p'}\textrm{ for all } f
\in L^{p'}(G)
\end{equation*}
and all conjugate exponents $p$ and $p'$ with $p' \in (1,2]$.
\end{proposition}
For a proof of this see, for example, Chapter 5 of Rudin
\cite{WR}. \item \emph{(Almost-orthogonality lemma)} To prove
Proposition \ref{local Bessels bound}, Green and Tao localized
Bessel's inequality to Bohr sets by using the following
almost-orthogonality lemma.
\begin{lemma}
\emph{(Cotlar's almost orthogonality lemma)} Suppose that $v$ and
$(w_j)$ are elements of a complex inner product space. Then
\begin{equation*}
\sum_j{|\langle v,w_j \rangle|^2} \leq \langle v,v\rangle
\max_j{\sum_i{|\langle w_i,w_j\rangle|}}.
\end{equation*}
\end{lemma} \item \emph{(Smoothed measures)} Suppose that $B(\Gamma,\delta)$ is a
regular Bohr set. We produce a range of smoothed alternatives to
the measure $\beta$; specifically suppose that $L \in \mathbb{N}$
and $\kappa \in (0,1]$. Then we may define
\begin{equation*}
\tilde{\beta}^{L,\kappa}_{\Gamma,\delta}:=\beta_{\Gamma,(1+\kappa)\delta}
\ast \beta_{\Gamma,\kappa\delta/L}^L,
\end{equation*}
where $\beta_{\Gamma,\kappa\delta/L}^L$ denotes the convolution of
$\beta_{\Gamma,\kappa\delta/L}$ with itself $L$ times. This
measure has the property that it is supported on
$B(\Gamma,(1+2\kappa)\delta)$ and uniform on $B(\Gamma,\delta)$,
indeed
\begin{equation}\label{equiv}
\tilde{\beta}_{\Gamma,\delta}^{L,\kappa}|_{B(\Gamma,\delta)} =
\frac{\mu_G|_{B(\Gamma,\delta)}}{\mu_G(B(\Gamma,(1+\kappa)\delta))}
=
\frac{\mu_G(B(\Gamma,\delta))}{\mu_G(B(\Gamma,(1+\kappa)\delta))}.\beta_{\Gamma,\delta}.
\end{equation}
It follows that every $f \in L^1(\beta_{\Gamma,\delta})$ has
$\wh{fd\beta_{\Gamma,\delta}}$ well approximated by
$\wh{fd\tilde{\beta}_{\Gamma,\delta}^{L,\kappa}}$. Specifically
\begin{equation}\label{estimator}
\wh{fd\tilde{\beta}_{\Gamma,\delta}^{L,\kappa}}(\gamma)=
(1+O(\kappa d))\wh{fd\beta_{\Gamma,\delta}}(\gamma)
\end{equation}
by regularity of $B(\Gamma,\delta)$.
\end{itemize}
We use almost-orthogonality and the smoothed measures to show the
following localization of Rudin's inequality. The proof of the
lemma to which this section is devoted then follows the usual
proof of Chang's theorem.
\begin{lemma}\label{local dual rudin} Suppose
that $G$ is a compact Abelian group and $B(\Gamma,\delta)$ is a
regular Bohr set. Suppose that $\Lambda$ is a set of characters.
Then there is a $\delta' \gg \delta/d|\Lambda|$ regular for
$\Gamma$ such that if $\Lambda$ is $\{\gamma:|\wh{\beta'}(\gamma)|
\geq 1/3\}$-dissociated then
\begin{equation*}
\|\wh{fd\beta}|_\Lambda\|_2 \ll \sqrt{p}\|f\|_{L^{p'}(\beta)}
\textrm{ for all } f \in L^{p'}(\beta)
\end{equation*}
and all conjugate exponents $p$ and $p'$ with $p' \in (1,2]$.
\end{lemma}
\begin{proof}
Begin by fixing the level of smoothing (i.e. the parameters $\kappa$
and $L$ of $\tilde{\beta}_{\Gamma,\delta}^{L,\kappa}$) that we
require and write $\tilde{\beta}$ for
$\tilde{\beta}_{\Gamma,\delta}^{L,\kappa}$. Set $L:=2k$ and recall
(\ref{estimator}):
\begin{equation*}
\wh{gd\tilde{\beta}}(\gamma)= (1+O(\kappa d))\wh{gd\beta}(\gamma)
\textrm{ for all } g \in L^1(\beta);
\end{equation*}
so we can pick $\kappa' \gg d^{-1}$ such that for all $\kappa \leq
\kappa'$
\begin{equation*} \frac{1}{2}|\wh{gd\beta}(\gamma)| \leq
|\wh{gd\tilde{\beta}}(\gamma)| \leq
\frac{3}{2}|\wh{gd\beta}(\gamma)| \textrm{ for all } g \in
L^1(\beta).
\end{equation*}
By Proposition \ref{ubreg}, we can take $\kappa$ with $\kappa'
\geq \kappa \gg d^{-1}$ such that $\delta':=\kappa\delta/L$ is
regular.

Define the Riesz product
\begin{equation*}
q(x):=\prod_{\lambda \in \Lambda}{\left(1+\frac{\lambda(x) +
\overline{\lambda}(x)}{2}\right)}.
\end{equation*}
Every term in the product is non-negative and so $q$ is
non-negative and it is fairly easy to compute the Fourier
transform of $q$:
\begin{equation}\label{ftq}
\wh{q}(\gamma)=\sum_{m:m.\Lambda=\gamma}{2^{-|m|}}.
\end{equation}
Since $\Lambda$ is $\{\gamma: |\wh{\beta'}(\gamma)| \geq
1/3\}$-dissociated, it is certainly vanilla dissociated and hence
$\wh{q}(0_{\wh{G}})=1$ and so, by non-negativity of $q$,
$\|q\|_1=1$.

Use $q$ to define the map
\begin{equation*}
T: L^1(\beta) \rightarrow L^1(G);g  \mapsto (gd\beta)\ast q,
\end{equation*}
and note that
\begin{equation*}
\|Tg\|_1 = \|(gd\beta) \ast q\|_1 \leq \|g\|_{L^1(\beta)}\|q\|_1 =
\|g\|_{L^1(\beta)}
\end{equation*}
by the triangle inequality. We now claim a corresponding result
for $\|Tg\|_2$, the proof of which we defer until we have finished
the proof of the lemma.
\begin{claim}
If $g \in L^2(\beta)$ then $\|Tg\|_2 \ll \|g\|_{L^2(\beta)}$.
\end{claim}
Assuming this claim, by the Riesz-Thorin interpolation theorem we
have
\begin{equation}\label{rtp}
\|Tg\|_{p'} \ll \|g\|_{L^{p'}(\beta)} \textrm{ for any }p' \in
[1,2].
\end{equation}
Hence, if $f \in L^{p'}(\beta)$,
\begin{eqnarray*}
\frac{1}{2}\|\wh{fd\beta}|_{\Lambda}\|_2 & \leq &
\|\wh{fd\beta}\wh{q}|_{\Lambda}\|_2  \textrm{ since }
\wh{q}(\lambda) \geq 1/2\textrm{ for all } \lambda \in \Lambda,\\
& = & \|\wh{Tf}|_\Lambda\|_2 \textrm{ by the definition of } T,\\
& \ll & \sqrt{p}\|Tf\|_{p'} \textrm{ by Rudin's inequality,}\\ &
\ll & \sqrt{p}\|f\|_{L^{p'}(\beta)} \textrm{ by (\ref{rtp})}.
\end{eqnarray*}
The lemma follows. It remains to prove the claim.
\begin{proof}[Proof of Claim.]
Begin by noting the following consequence of (\ref{equiv}).
\begin{equation}\label{str}
\|Tg\|_2^2=\left(\frac{\mu_G(B(\Gamma,\delta(1+\kappa)))}{\mu_G(B(\Gamma,\delta))}\right)^2\|(gd\tilde{\beta})
\ast q\|_2^2.
\end{equation}
By Plancherel's theorem
\begin{equation*}
\|(gd\tilde{\beta}) \ast q\|^2 = \sum_{\gamma \in
\wh{G}}{|\wh{(gd\tilde{\beta})}(\gamma)\wh{q}(\gamma)|^2} =
\sum_{\gamma \in \wh{G}}{|\langle
g,\wh{q}(\gamma)\gamma\rangle_{L^2(\tilde{\beta})}|^2}.
\end{equation*}
Cotlar's Almost Orthoginality Lemma applied to the second sum
gives
\begin{eqnarray*}
\|(gd\tilde{\beta}) \ast q\|^2 & \leq & \langle
g,g\rangle_{L^2(d\tilde{\beta})}\max_{\gamma}{\sum_{\gamma'}{|\langle
\wh{q}(\gamma)\gamma,\wh{q}(\gamma')\gamma'\rangle_{L^2(\tilde{\beta})}|}}\\
& \leq &
\|g\|_{L^2(d\tilde{\beta})}^2\max_{\gamma}{\sum_{\gamma'}{\wh{q}(\gamma')|\wh{\tilde{\beta}}(\gamma-\gamma')|}}.
\end{eqnarray*}
For any $\gamma \in \wh{G}$ we can estimate the last sum in a
manner independent of $\gamma$ by using a positivity argument:
\begin{eqnarray*}
\sum_{\gamma' \in
\wh{G}}{\wh{q}(\gamma')|\wh{\tilde{\beta}}(\gamma - \gamma')|}  &
= & \sum_{\gamma' \in
\wh{G}}{\wh{q}(\gamma-\gamma')|\wh{\beta}(\gamma')
\wh{\beta'}(\gamma')^L|} \textrm{ by definition of
$\tilde{\beta}$,}\\ & \leq & \sum_{\gamma' \in \wh{G}
}{\wh{q}(\gamma-\gamma')|\wh{\beta'}(\gamma')|^L} \textrm{ since
$|\wh{\beta}(\gamma')| \leq \|\beta\| = 1$ and $\wh{q}\geq 0$, }\\
& = & \wh{qd\beta'^L}(\gamma) \textrm{ since $L$ is even and
$\wh{q} \geq 0$,}\\ & \leq & \|q\|_{L^1(\beta'^L)} \\ & =
&\wh{qd\beta'^L}(0_{\wh{G}}) \textrm{ by non-negativity of
$qd\beta'^L$,}
\\  & = & \sum_{\gamma' \in
\wh{G}}{\wh{q}(\gamma')|\wh{\beta'}(\gamma')|^L} \textrm{ by
symmetry of $\wh{q}$.}
\end{eqnarray*}
We estimate this in turn by splitting the range of summation into
two parts:
\begin{equation}\label{crty} \sum_{\gamma' \in
\wh{G}}{\wh{q}(\gamma')|\wh{\beta'}(\gamma')|^L} \leq
\sum_{\gamma':|\wh{\beta'}(\gamma')| \geq
1/3}{\wh{q}(\gamma')|\wh{\beta'}(\gamma')|^L} +
\sum_{\gamma':|\wh{\beta'}(\lambda)| \leq
1/3}{\wh{q}(\gamma')|\wh{\beta'}(\gamma')|^L}.
\end{equation}
\begin{enumerate}
\item For the first sum: $|\wh{q}(\gamma')| \leq \|q\|_1 = 1$ and
$|\wh{\beta'}(\gamma')^L| \leq \|\beta'^L\| = 1$ so that each
summand is at most 1, furthermore $\supp \wh{q} \subset \langle
\Lambda \rangle$ so
\begin{equation*}
\sum_{\gamma':|\wh{\beta'}(\gamma')| \geq
1/3}{\wh{q}(\gamma')|\wh{\beta'}(\gamma')|^L} \leq \sum_{\gamma'
\in \langle \Lambda \rangle:|\wh{\beta'}(\gamma')| \geq 1/3}{1}.
\end{equation*}
This range of summation contains at most 1 element by
$\{\gamma:|\wh{\beta'}(\gamma)|\geq 1/3\}$-dissociativity of
$\Lambda$, and hence the sum is bounded above by 1. \item For the
second sum: $|\wh{q}(\gamma')| \leq \|q\|_1 =1$ and
$|\wh{\beta'}(\gamma')^L| \leq 3^{-L}$ for $\gamma'$ in the range
of summation so that each summand is at most $9^{-|\Lambda|}$,
however $\supp \wh{q} \subset \langle \Lambda \rangle$ and
$|\langle \Lambda \rangle| \leq 3^{|\Lambda|}$ so
\begin{equation*}
\sum_{\gamma':|\wh{\beta'}(\gamma')| \leq
1/3}{\wh{q}(\gamma')|\wh{\beta'}(\gamma')|^L} \leq \sum_{\gamma'
\in \langle \Lambda \rangle}{9^{-|\Lambda|}} \leq 1.
\end{equation*}
\end{enumerate}
It follows that the right hand side of (\ref{crty}) is bounded
above by 2 and hence that
\begin{equation*}
\|(gd\tilde{\beta}) \ast q\|^2 \leq
2\|g\|_{L^2(d\tilde{\beta})}^2.
\end{equation*}
This, (\ref{str}) and (\ref{equiv}) yield
\begin{equation*}
\|Tg\|_2^2 \leq
2\frac{\mu_G(B(\Gamma,\delta(1+\kappa)))}{\mu_G(B(\Gamma,\delta))}\|g\|_{L^2(\beta)}^2,
\end{equation*}
from which the result follows by regularity.
\end{proof}
\end{proof}

\begin{proof}[Proof of Lemma \ref{content chang}.] By the localized dual of Rudin's inequality (Lemma \ref{local dual
rudin}) for any $p' \in (1,2]$ we have
\begin{equation*}
|\Lambda|.\epsilon^2\|f\|_{L^1(\beta)}^2 \leq \sum_{\lambda \in
\Lambda}{|\wh{fd\beta}(\lambda)|^2} =
\|\wh{fd\beta}|_{\Lambda}\|_2^2 \ll p\|f\|_{L^{p'}(\beta)}^2,
\end{equation*}
where $p$ is the conjugate exponent of $p'$. The $\log$-convexity
of $\|.\|_{L^{p'}(\beta)}$ gives
\begin{equation*}
|\Lambda| \ll
\epsilon^{-2}p\left(\frac{\|f\|_{L^2(\beta)}}{\|f\|_{L^1(\beta)}}\right)^{\frac{4}{p}}.
\end{equation*}
Optimizing $p$ gives the result.
\end{proof}

\section{Local Fourier analysis and the iteration method}

The tools of local Fourier analysis were originally developed with
iteration in mind. Specifically if $A$ is a subset of a regular
Bohr set $B$ then we shall often have an argument which tells us
that there is a large $\ell^p$-mass of the local Fourier transform
$\wh{\chi_Ad\beta}$ and, as is the case in the non-local setting,
this leads to a density increment on a (sub-)Bohr neighborhood.
For our purposes we have the following two standard lemmas which
take a large $\ell^\infty$ and $\ell^2$ Fourier-space mass and
convert it into a density increment.

\begin{lemma}\label{linfty density}
\emph{($\ell^\infty$ density increment argument)} Suppose that $G$
is a compact Abelian group and $B=B(\Gamma,\delta)$ is a regular
Bohr set. Suppose that $A \subset B$ has relative density $\alpha$
and write $f:=\chi_A - \alpha \chi_B$. Suppose that
\begin{equation*}
|\wh{fd\beta}(\gamma)| \geq \eta \alpha \textrm{ for some } \gamma
\in \wh{G}.
\end{equation*}
Then there is a regular Bohr set $B':=B(\Gamma',\delta')$ with
$\Gamma':=\Gamma\cup\{\gamma\}$ and $\delta \geq \delta' \gg
\eta\alpha\delta/d$ such that
\begin{equation*}
\|\chi_A \ast \beta'\|_{L^\infty(\beta)} \geq
\alpha\left(1+2^{-3}\eta\right).
\end{equation*}
\end{lemma}
\begin{proof} Let $\delta' \in (0,1]$ be a constant to be determined later. A trivial instance of Hausdorff's
inequality tells us that
\begin{equation}\label{y}
\|(fd\beta) \ast \beta'\| \geq
|\wh{fd\beta}(\gamma)||\wh{\beta'}(\gamma)| \geq
\eta\alpha|\wh{\beta'}(\gamma)|.
\end{equation}
Since $B' \subset B(\{\gamma\},\delta')$ we have
$|\wh{\beta'}(\gamma)| \geq 1-O(\delta')$. It follows that there
is a $\delta'_0 \gg 1$ such that if $\delta' \leq \delta'_0$ then
$|\wh{\beta'}(\gamma)| \geq 1/2$. Now
\begin{equation*}
\int{d((fd\beta) \ast \beta')}=0,
\end{equation*}
hence by (\ref{y})
\begin{equation*}
\int{d((fd\beta) \ast \beta')_+} \geq
\eta\alpha|\wh{\beta'}(\gamma)|/2 \geq \eta\alpha/4.
\end{equation*}
If follows from the regularity of $B$ and the fact that $B'
\subset B(\Gamma,\delta')$ that
\begin{equation*}
\| (fd\beta) \ast \beta' - (f\ast \beta') d\beta\| =
O(d\delta'\delta^{-1}),
\end{equation*}
and so
\begin{equation*}
\int{(f\ast \beta')_+d\beta} \geq \eta\alpha/4 +
O(d\delta'\delta^{-1}).
\end{equation*}
By regularity of $B(\Gamma,\delta)$ we have
\begin{equation*}
\int{(f\ast \beta')_+d\beta} \leq \|\chi_A \ast \beta'\|_\infty -
\alpha + O(d\delta'\delta^{-1}),
\end{equation*}
so
\begin{equation*}
\|\chi_A \ast \beta'\|_\infty \geq \alpha(1+1/4) +
O(d\delta'\delta^{-1}).
\end{equation*}
Hence by Proposition \ref{ubreg} we can pick $\delta' \gg
\eta\alpha\delta/d$ regular for $\Gamma'$ with $\delta' \leq
\delta'_0$ and $\delta' \leq \delta$, and such that the conclusion
of the lemma holds.
\end{proof}

\begin{lemma}
\label{l2 density} \emph{($\ell^2$ density increment argument)}
Suppose that $G$ is a compact Abelian group and
$B=B(\Gamma,\delta)$ is a regular Bohr set. Suppose that $A
\subset B$ has relative density $\alpha$ and write $f:=\chi_A -
\alpha \chi_B$. Suppose that $B'=B(\Gamma',\delta')$ is a Bohr set
with $\Gamma \subset \Gamma'$ and
\begin{equation*}
\langle f \ast \beta',(fd\beta) \ast \beta' \rangle \geq
c\alpha^2.
\end{equation*}
Then
\begin{equation*}
\|\chi_A \ast \beta'\|_\infty \geq \alpha(1+c) +
O(d\delta'\delta^{-1}).
\end{equation*}
\end{lemma}
\begin{proof}
We expand the inner product:
\begin{eqnarray*}
\langle f \ast \beta',(fd\beta) \ast \beta' \rangle & = & \langle
\chi_A \ast \beta', (\chi_Ad\beta) \ast \beta' \rangle -
\alpha\langle \chi_B \ast \beta', (\chi_Ad\beta)\ast \beta'\rangle
\\ & & - \alpha \langle \chi_A \ast \beta', \beta \ast \beta'
\rangle + \alpha^2\langle \chi_B \ast \beta', \beta \ast \beta'
\rangle.
\end{eqnarray*}
Now we estimate each term. First
\begin{eqnarray*}
\langle \chi_A \ast \beta', (\chi_Ad\beta) \ast \beta' \rangle &
\leq & \|\chi_A \ast \beta'\|_\infty \|(\chi_Ad\beta) \ast
\beta'\|\\ & \leq &\|\chi_A \ast
\beta'\|_\infty\|\chi_A\|_{L^1(\beta)}\|\beta'\|= \|\chi_A \ast
\beta'\|_\infty\alpha.
\end{eqnarray*}
By Lemma \ref{contlem} and the fact that $B(\Gamma',\delta')
\subset B(\Gamma,\delta')$ we have
\begin{equation*}
\|\beta \ast \beta' \ast \beta' - \beta\| = O(d\delta'\delta^{-1})
\end{equation*}
whence
\begin{equation*}
\langle \chi_B \ast \beta', (\chi_Ad\beta)\ast \beta'\rangle =
\langle \beta \ast \beta' \ast \beta', \chi_A \rangle = \alpha +
O(d\delta'\delta^{-1}),
\end{equation*}
\begin{equation*}
\langle \chi_A \ast \beta', \beta \ast \beta' \rangle = \langle
\chi_A, \beta \ast \beta'\ast \beta' \rangle = \alpha +
O(d\delta'\delta^{-1}),
\end{equation*}
and
\begin{equation*}
\langle \chi_B \ast \beta', \beta \ast \beta' \rangle = \langle
\chi_B , \beta \ast \beta' \ast \beta' \rangle = 1 +
O(d\delta'\delta^{-1}).
\end{equation*}
It follows that
\begin{equation*}
\alpha\|\chi_A \ast \beta'\|_\infty \geq \alpha^2(1+c) + O(\alpha
d\delta'\delta^{-1}),
\end{equation*}
from which we get the result on division by $\alpha$.
\end{proof}

\section{Proof of Theorem \ref{green1new}}

We begin with a brief overview of the argument in the model
setting. This argument can be made to prove the following result,
which was first established by Green in \cite{BJGRKP}.
\begin{theorem}\label{modelaa}
Suppose that $G$ is a finite dimensional compact vector space over
$\mathbb{F}_2$ and $A \subset G$ has density $\alpha>0$. Then
$A+A$ contains a subspace of dimension $2^{-4}\alpha^2\dim G$.
\end{theorem}
There are three main ingredients to the proof of this result.
First we have the iteration lemma - the driving force. In words it
says that either $A+A$ contains most of $G$ or we can find an
affine subspace on which $A$ has increased density.
\begin{lemma}
\emph{(Model iteration lemma)} Suppose that $G$ is a compact vector
space over $\mathbb{F}_2$. Suppose that $A \subset G$ has density
$\alpha$. Suppose that $\sigma \in (0,1]$. Then at least one of the
following is true.
\begin{enumerate}
\item (The sumset contains most of $G$) $A+A$ contains at least a
proportion $1-\sigma$ of $G$. \item (Density increment) There is a
subspace $V$ of $G$ such that
\begin{equation*}
\|\chi_{A} \ast \mu_V\|_\infty \geq \alpha(1+1/4) \textrm{ and }
\cod V \leq 8\alpha^{-2}\log \sigma^{-1}.
\end{equation*}
\end{enumerate}
\end{lemma}
The proof of this is not difficult; we sketch the main ideas now.
It was a crucial insight of Green in \cite{BJGAA} to get control
of $A \subset G$ by looking at its complement. Specifically if $S
\subset (A+A)^c$ then we have
\begin{equation*}
\langle \chi_A \ast \chi_A, \chi_S \rangle=0.
\end{equation*}
Green employed an ingenious argument to exploit this information;
ours is less sophisticated. Plancherel's theorem and the triangle
inequality in the usual fashion will give
\begin{equation*}
\sum_{\gamma \neq
0_{\wh{G}}}{|\wh{\chi_A}(\gamma)|^2|\wh{\chi_S}(\gamma)|} \geq
\alpha^2 \sigma,
\end{equation*}
so if $\mathcal{L}$ is the set of non-trivial characters
supporting large values of $|\wh{\chi_S}|$ then it follows easily
enough that
\begin{equation*}
\sum_{\gamma \in \mathcal{L}}{|\wh{\chi_A}(\gamma)|^2} \gg
\alpha^2.
\end{equation*}
Such a bound provides an $\ell^2$ density increment for $A$; we
bound the codimension of the subspace on which we get the
increment by using Chang's theorem.

The second ingredient is a simple pigeonhole argument which says
that if a set contains a large proportion of a vector space then it
must contain a large affine subspace.
\begin{lemma}
\emph{(Pigeonhole lemma)} Suppose that $G$ is a finite dimensional
compact vector space over $\mathbb{F}_2$ and that $A \subset G$ has
density $\alpha>1-\sigma$. Then $A$ contains a coset of a subspace
of dimension $\lfloor\log_2 \sigma^{-1}\rfloor$ provided $G$
contains a subspace of dimension $\lfloor\log_2 \sigma^{-1}\rfloor$.
\end{lemma}

The iteration necessary to prove Theorem \ref{modelaa} is now very
simple. At each stage of the argument we apply the iteration lemma
and conclude that either $A+A$ contains a large portion of an affine
space or the density of $A$ can be increased on an affine subspace.
The density of $A$ cannot be increased forever, so eventually $A+A$
contains a large portion of an affine space and so we may apply the
pigeonhole lemma to conclude that $A+A$ contains a large affine
space. Optimizing the parameter $\sigma$ gives the result.

We now turn to the matter of transferring these ideas to the
general setting.
\begin{lemma}\label{itlem}\emph{(Iteration lemma)}
Suppose that $G$ is a compact Abelian group and $B(\Gamma,\delta)$
is a regular Bohr set. Suppose that $A_1,A_2 \subset
B(\Gamma,\delta)$. Write $\alpha$ for the geometric mean of the
densities of $A_1$ and $A_2$ in $B(\Gamma,\delta)$. Suppose that
$\sigma \in (0,1]$. Then at least one of the following is true.
\begin{enumerate}
\item (The sumset contains most of a Bohr set) There is a regular
Bohr set $B(\Gamma,\delta')$ such that $A_1+A_2$ contains at least
a proportion $1-\sigma$ of $B(\Gamma,\delta')$ and $\delta' \gg
\alpha^4\delta/d$. \item (Density increment) There is a regular
Bohr set $B(\Gamma\cup\Lambda,\delta'')$ such that
\begin{equation*}
\|\chi_{A_1} \ast
\beta_{\Gamma\cup\Lambda,\delta''}\|_\infty\|\chi_{A_2}\ast
\beta_{\Gamma\cup\Lambda,\delta''}\|_\infty \geq
\alpha^2(1+2^{-4}),
\end{equation*}
and
\begin{equation*}
|\Lambda| \ll \alpha^{-2}\log \sigma^{-1} \textrm{ and }\delta''
\gg \delta\alpha^6/|\Gamma|^3\log \sigma^{-1}.
\end{equation*}
\end{enumerate}
\end{lemma}
As is typical of arguments of this type the proof is quite
technical; to simplify the presentation we extract two lemmas from
the main argument and place them at the end. Pedagogically, it would
be most appropriate to present them now, but they are hard to
motivate without following main proof; hence the order we have
chosen.
\begin{proof}
We write $\alpha_1$ and $\alpha_2$ for the densities of $A_1$ and
$A_2$ (respectively) in $B(\Gamma,\delta)$, and $d$ for the size of
$\Gamma$. We may assume that $\alpha_1,\alpha_2>0$ since otherwise
the result is trivial.

Let $\delta'\in (0,1]$ be a constant, the value of which will fall
out of the proof and write $B'$ for the Bohr set $B(\Gamma,\delta')$
and $B$ for the Bohr set $B(\Gamma,\delta)$. Either we are in the
first case of the lemma or we may pick $S \subset B' \setminus
(A_1+A_2)$ with $\beta'(S)=\sigma$. We have
\begin{equation}\label{zerohere}
\langle \chi_{A_1} \ast (\chi_{A_2} d\beta),
\chi_S\rangle_{L^2(\beta')}=0.
\end{equation}
Write $f_i$ for the \emph{balanced function} $\chi_{A_i} -
\alpha_i\chi_B$ of $A_i$ in $B$. Then
\begin{equation*}
f_1 \ast (f_2d\beta) = \chi_{A_1} \ast (\chi_{A_2}d\beta) -
\alpha_1\chi_{B} \ast (\chi_{A_2}d\beta) -\chi_{A_1} \ast
\alpha_2\beta + \alpha_1\alpha_2\chi_B \ast \beta.
\end{equation*}
For $x \in B'$, the last three terms on the right may be estimated
using Corollary \ref{contlemcor2}:
\begin{equation*}
\alpha_1\chi_B \ast (\chi_{A_2}d\beta)(x) = \alpha_1\chi_{A_2}\ast
\beta(x)  = \alpha_1(\alpha_2+O(d\delta'\delta^{-1}));
\end{equation*}
\begin{equation*}
\chi_{A_1} \ast \alpha_2\beta (x)= \alpha_2(\alpha_1 +
O(d\delta'\delta^{-1}));
\end{equation*}
\begin{equation*}
\alpha_1\alpha_2 \chi_B \ast
\beta(x)=\alpha_1\alpha_2(1+O(d\delta'\delta^{-1})).
\end{equation*}
Whence
\begin{equation*}
f_1 \ast (f_2d\beta) = \chi_{A_1} \ast (\chi_{A_2}d\beta)-
\alpha^2+O(d\delta'\delta^{-1}).
\end{equation*}
It follows from this and (\ref{zerohere}) that
\begin{equation*}
\langle f_1 \ast (f_2d\beta), \chi_S \rangle_{L^2(\beta')}= -
\alpha^2\sigma+O(d\delta'\delta^{-1}\sigma).
\end{equation*}
Apply Plancherel's theorem to this inner product to produce a
Fourier statement:
\begin{equation*}
\sum_{\gamma \in
\wh{G}}{\wh{f_1}(\gamma)\wh{f_2d\beta}(\gamma)\overline{\wh{\chi_Sd\beta'}(\gamma)}}=-\alpha^2\sigma+O(d\delta'\delta^{-1}\sigma),
\end{equation*}
so by the triangle inequality
\begin{equation}\label{triangle}
\sum_{\gamma \in
\wh{G}}{|\wh{f_1}(\gamma)||\wh{f_2d\beta}(\gamma)||\wh{\chi_Sd\beta'}(\gamma)|}\geq\alpha^2\sigma+O(d\delta'\delta^{-1}\sigma).
\end{equation}

Let $\mathcal{L}$ be the set of characters supporting the large
values of $|\wh{\chi_Sd\beta'}|$:
\begin{equation*}
\mathcal{L}:=\{\gamma \in\wh{G}:|\wh{\chi_Sd\beta'}(\gamma)| \geq
\alpha\sigma/2\}.
\end{equation*}
The characters supporting small values of $|\wh{\chi_Sd\beta'}|$
only support a small amount of the sum in (\ref{triangle});
specifically by Lemma \ref{techlem1} applied with
$h=|\wh{\chi_Sd\beta'}(\gamma)|\chi_{\mathcal{L}^c}(\gamma)$ we have
\begin{equation*}
\sum_{\gamma \not \in
\mathcal{L}}{|\wh{f_1}(\gamma)||\wh{f_2d\beta}(\gamma)||\wh{\chi_Sd\beta'}(\gamma)|}<
\frac{\alpha^2\sigma}{2}.
\end{equation*}
Inserting this into (\ref{triangle}) we conclude that
\begin{equation*}
\sum_{\gamma \in
\mathcal{L}}{|\wh{f_1}(\gamma)||\wh{f_2d\beta}(\gamma)||\wh{\chi_Sd\beta'}(\gamma)|}\geq
\frac{\alpha^2\sigma}{2}+O(d\delta'\delta^{-1}\sigma).
\end{equation*}
We have the trivial inequality $|\wh{\chi_Sd\beta'}(\gamma)| \leq
\sigma$ and so (since $\sigma>0$) we divide through by $\sigma$ to
get
\begin{equation*}
\sum_{\gamma \in
\mathcal{L}}{|\wh{f_1}(\gamma)||\wh{f_2d\beta}(\gamma)|} \geq
\frac{\alpha^2}{2}+O(d\delta'\delta^{-1}).
\end{equation*}

By Proposition \ref{local Changs bound} there is a set of
characters $\Lambda$ and a $\delta''_0 \in (0,1]$ with
\begin{equation*}
|\Lambda| \ll \alpha^{-2} \log \sigma^{-1} \textrm{ and }
\delta''_0 \gg \frac{\delta'\alpha^2}{d^2 \log \sigma^{-1}}
\end{equation*}
such that
\begin{equation*}
\mathcal{L} \subset \{\gamma:|1-\gamma(x)| \leq 1/2 \textrm{ for
all } x \in B(\Gamma\cup \Lambda,\delta''_0)\}.
\end{equation*}
Suppose that $\delta'' \leq \delta''_0$, and write $\beta''$ for
$\beta_{\Gamma\cup\Lambda,\delta''}$. If $\gamma \in \mathcal{L}$
then $|\wh{\beta''}(\gamma)|\geq 1/2$, so
\begin{equation}\label{near last estimate}
\sum_{\gamma \in
\wh{G}}{|\wh{f_1}(\gamma)||\wh{f_2d\beta}(\gamma)||\wh{\beta''}(\gamma)|^2}\geq
2^{-3}\alpha^2+O(d\delta'\delta^{-1}).
\end{equation}
Now, by Lemma \ref{techlem2} we have
\begin{equation*}
\sum_{\gamma \in
\wh{G}}{|\wh{f_1}(\gamma)||\wh{f_2d\beta}(\gamma)||\wh{\beta''}(\gamma)|^2}
\leq \alpha^2 \max_{1 \leq i\leq 2}{\alpha_i^{-2}\langle f_i \ast
\beta'',(f_id\beta) \ast \beta''\rangle}.
\end{equation*}
Which, combined with (\ref{near last estimate}), ensures that there
is some $k$ with $1 \leq k \leq 2$ such that
\begin{equation*}
\langle f_k \ast \beta'',(f_kd\beta) \ast \beta''\rangle \geq
\alpha_k^2\left(2^{-3}+O(d\delta'\delta^{-1}\alpha^{-2})\right).
\end{equation*}
Now we can apply Lemma \ref{l2 density} to get that
\begin{equation*}
\|\chi_{A_k} \ast \beta''\|_\infty \geq \alpha_k\left(1+2^{-3} +
O(d\delta'\delta^{-1}\alpha^{-2})\right) +
O(d\delta''\delta^{-1}).
\end{equation*}
However for $1 \leq i \leq 2$ we have
\begin{eqnarray*}
\|\chi_{A_i} \ast \beta''\|_\infty & \geq &\|\chi_{A_i} \ast
\beta''\|_{L^1(\beta)}\\ & = & \int{\chi_{A_i} d(\beta \ast
\beta'')}\\ &  = & \alpha_i + O(d\delta''\delta^{-1}) \textrm{ by
Corollary \ref{contlemcor} since } \supp \beta'' \subset
B(\Gamma,\delta''),
\end{eqnarray*}
so that
\begin{equation*}
\|\chi_{A_1} \ast \beta''\|_\infty\|\chi_{A_2} \ast
\beta''\|_\infty \geq \alpha^2(1+2^{-3}) +
O(d\delta'\delta^{-1}\alpha^{-2})+ O(d\delta''\delta^{-1}).
\end{equation*}
Assume that $\delta'' \leq \delta'$, so that
\begin{equation*}
\|\chi_{A_1} \ast \beta''\|_\infty\|\chi_{A_2} \ast
\beta''\|_\infty \geq \alpha^2(1+2^{-3}) +
O(d\delta'\delta^{-1}\alpha^{-2}).
\end{equation*}
We now pick $\delta'$ regular for $\Gamma$ such that the error
term in the above expression is at most $2^{-4}\alpha^2$. This can
be done by Proposition \ref{ubreg} whilst keeping $\delta' \gg
\alpha^4 \delta/d$. Finally we pick $\delta''$ regular for
$\Gamma\cup\Lambda$ subject to the two assumptions of $\delta''
\leq \delta'$ and $\delta'' \leq \delta''_0$. This can be done by
Proposition \ref{ubreg} whilst keeping $\delta'' \gg
\alpha^6\delta/d^3\log \sigma^{-1}$, and so the lemma is proved.
\end{proof}
We now prove the two technical claims we required.
\begin{lemma}\label{techlem1}
  Suppose that $G$ is a compact Abelian group and $B(\Gamma,\delta)$
is a Bohr set. Suppose that $f_1,f_2 \in L^2(\beta)$ and $h \in
\ell^\infty(\wh{G})$. Then
\begin{equation*}
  \sum_{\gamma \in
  \wh{G}}{|\wh{f_1}(\gamma)||\wh{f_2d\beta}(\gamma)|h(\gamma)}
  \leq
  \|f_1\|_{L^2(\beta)}\|f_2\|_{L^2(\beta)}\|h\|_{\ell^\infty(\wh{G})}.
\end{equation*}
\end{lemma}
\begin{proof}
Start with the fact that
\begin{equation*}
\sum_{\gamma \in
\wh{G}}{|\wh{f_1}(\gamma)||\wh{f_2d\beta}(\gamma)|h(\gamma)}
 < \|h\|_{\ell^\infty(\wh{G})}\sum_{\gamma \in
\wh{G}}{|\wh{f_1}(\gamma)||\wh{f_2d\beta}(\gamma)|}.
\end{equation*}
We estimate the sum on the right using the Cauchy-Schwarz inequality
and Plancherel's theorem.
\begin{eqnarray*}
\sum_{\gamma \in \wh{G}}{|\wh{f_1}(\gamma)||\wh{f_2d\beta}(\gamma)|}
& \leq & \left(\sum_{\gamma \in
\wh{G}}{|\wh{f_1}(\gamma)|^2}\right)^{\frac{1}{2}}\left(\sum_{\gamma
\in \wh{G}}{|\wh{f_2d\beta}(\gamma)|^2}\right)^{\frac{1}{2}}\\
& = & \|f_1\|_2\|f_2d\beta\|_2\\ & = &
\|f_1\|_{L^2(\beta)}\|f_2\|_{L^2(\beta)} \textrm{ since $\beta$ is
uniform on $B(\Gamma,\delta)$.}
\end{eqnarray*}
Putting these two inequalities together gives the result.
\end{proof}
Similarly we have the following.
\begin{lemma}\label{techlem2}
  Suppose that $G$ is a compact Abelian group and $B(\Gamma,\delta)$
  and $B(\Gamma',\delta')$ are Bohr sets. Suppose that $f_1,f_2 \in L^2(\beta)$. Then
\begin{equation*}
  \sum_{\gamma \in \wh{G}}{|\wh{f_1}(\gamma)||\wh{f_2d\beta}(\gamma)||\beta'(\gamma)|^2}
  \leq \|f_1\|_{L^2(\beta)}\|f_2\|_{L^2(\beta)}\max_{1 \leq i \leq
2}{\|f_i\|_{L^2(\beta)}^{-2}\langle f_i \ast \beta',(f_id\beta) \ast
\beta' \rangle}.
\end{equation*}
\end{lemma}
\begin{proof}
Apply the Cauchy-Schwarz inequality to the sum on the left to bound
it above by
\begin{equation*}
\left(\sum_{\gamma \in
\wh{G}}{|\wh{f_1}(\gamma)|^2|\wh{\beta'}(\gamma)|^2}\right)^\frac{1}{2}\left(\sum_{\gamma
\in
\wh{G}}{|\wh{f_2d\beta}(\gamma)|^2|\wh{\beta'}(\gamma)|^2}\right)^{\frac{1}{2}};
\end{equation*}
as before this can be rewritten as
\begin{equation*}
\left( \sum_{\gamma \in
\wh{G}}{\wh{f_1}(\gamma)\wh{\beta''}(\gamma)\overline{\wh{f_1d\beta}(\gamma)\wh{\beta'}(\gamma)}}\right)^{\frac{1}{2}}\left(
\sum_{\gamma \in
\wh{G}}{\wh{f_2}(\gamma)\wh{\beta''}(\gamma)\overline{\wh{f_2d\beta}(\gamma)\wh{\beta'}(\gamma)}}\right)^{\frac{1}{2}}.
\end{equation*}
Now apply Plancherel's theorem to this to conclude that it is equal
to
\begin{equation*}
\langle f_1 \ast \beta',(f_1d\beta) \ast
\beta'\rangle^{\frac{1}{2}}\langle f_2 \ast \beta',(f_2d\beta) \ast
\beta'\rangle^{\frac{1}{2}},
\end{equation*}
from which the lemma follows.
\end{proof}

The next lemma is a local version of the following easy application
of the pigeonhole principle: If $A \subset \mathbb{Z}/N\mathbb{Z}$
has density at least $1-\sigma$ then $A$ contains an arithmetic
progression of length roughly $\sigma^{-1}$. It turns out not to be
hard to localize this observation.

\begin{lemma}\label{phole}
Suppose that $G=\mathbb{Z}/N\mathbb{Z}$ and $B(\Gamma,\delta)$ is
a regular Bohr set. Suppose that $A \subset G$. Suppose that
$\sigma \in (0,1]$. Suppose that $A$ contains at least a
proportion $1-\sigma$ of $B(\Gamma,\delta)$. Then either
$\sigma^{-1} \gg d^{-1}\delta N^{\frac{1}{d}}$ or $A$ contains an
arithmetic progression of length at least $(4\sigma)^{-1}$.
\end{lemma}
\begin{proof}
Let $\eta$ be a constant to be optimized later.

First we find a large number of long arithmetic progressions in
$B(\Gamma,\delta)$, all with the same common difference. Pick $y
\neq 0$ from $B(\Gamma,2^{\frac{1}{d}}N^{-\frac{1}{d}})$; such a
$y$ certainly exists by Lemma \ref{bohrsize} which ensures that
$|B(\Gamma,2^{\frac{1}{d}}N^{-\frac{1}{d}})| \geq 2$. It follows
that
\begin{equation*}
x \in B(\Gamma,\delta(1-\eta)) \Rightarrow x,x+y,x+2y,...,x+Ly \in
B(\Gamma,\delta)
\end{equation*}
for $L \leq \eta\delta N^{\frac{1}{d}}2^{-\frac{1}{d}}$. Hence if
$(4\sigma)^{-1} \leq \eta\delta N^{\frac{1}{d}}2^{-\frac{1}{d}}$
then there are at least $\mu_G(B(\Gamma,\delta(1-\eta)))N$
arithmetic progressions of common difference $y$ and length
$(4\sigma)^{-1}$ in $B(\Gamma,\delta)$. Moreover, since the common
difference is the same for each progression, each point is in at
most $2.(4\sigma)^{-1}=(2\sigma)^{-1}$ of these progressions.

If $A$ does not contain any of these progressions then it misses
at least one point in each progression and hence at least
$\mu_G(B(\Gamma,\delta(1-\eta)))N/ (2\sigma)^{-1}$ points of
$B(\Gamma,\delta)$. It follows that
\begin{equation*}
1-\sigma \leq \int{\chi_A d\beta} \leq 1-
2\sigma\frac{\mu_G(B(\Gamma,\delta(1-\eta)))}{\mu_G(B(\Gamma,\delta))}.
\end{equation*}
By regularity of $\delta$ we can pick $\eta \gg d^{-1}$ such that
\begin{equation*}
\frac{\mu_G(B(\Gamma,\delta(1-\eta)))}{\mu_G(B(\Gamma,\delta))}
\geq \frac{2}{3},
\end{equation*}
from which we conclude that $1-\sigma \leq 1-4\sigma/3$, this
contradicts the fact that $\sigma$ is positive and so the lemma is
proved.
\end{proof}

Finally we put the previous two results together to prove Theorem
\ref{green1new}.

\begin{proof}
[Proof of Theorem \ref{green1new}.] Let $\sigma>0$ be a constant to
be optimized later. We construct a sequence of regular Bohr sets
$B(\Gamma_k,\delta_k)$ iteratively. Write
\begin{equation*}
\beta_k:=\beta_{\Gamma_k,\delta_k}, d_k:=|\Gamma_k| \textrm{ and }
\alpha_k:=\sqrt{\|\chi_{A_1} \ast \beta_{k}\|_\infty \|\chi_{A_2}
\ast \beta_{k}\|_\infty}.
\end{equation*}
We initialize the iteration with $\Gamma_0=\{0_{\wh{G}}\}$ and
$\delta_0 \gg 1$ regular for $\Gamma_0$ by Proposition
\ref{ubreg}.

Suppose that we are at stage $k$ of the iteration. Since
$B(\Gamma_k,\delta_k)$ has positive measure $\chi_{A_i} \ast
\beta_k$ is continuous and hence we make take $x_1$ and $x_2$ such
that
\begin{equation*}
\chi_{A_i} \ast \beta_k(x_i) = \|\chi_{A_i} \ast \beta_k\|_\infty
\textrm{ for } i=1,2.
\end{equation*}
Apply Lemma \ref{itlem} to the sets $(A_1-x_1)\cap
B(\Gamma_k,\delta_k)$ and $(A_2-x_2) \cap B(\Gamma_k,\delta_k)$
and the regular Bohr set $B(\Gamma_k,\delta_k)$.
\begin{enumerate}
\item Either $(A_1-x_1)+(A_2-x_2)$ contains at least a proportion
$1-\sigma$ of some regular Bohr set $B(\Gamma_k,\delta_k')$ with
$\delta_k' \gg \alpha_k^4\delta_k/d_k$. In which case we apply
Lemma \ref{phole} to conclude that either $\sigma^{-1} \gg
d_k^{-1}\delta_k'N^{\frac{1}{d_k}}$ or $A_1+A_2-x_1-x_2$ (and
hence $A_1+A_2$) contains an arithmetic progression of length
$(4\sigma)^{-1}$. \item Or there is a regular Bohr set
$B(\Gamma_{k+1},\delta_{k+1})$ such that
\begin{equation*}
\alpha_{k+1}^2 \geq \alpha_k^2(1+2^{-4}), \delta_{k+1} \gg
\frac{\alpha_k^6\delta_k}{d_k^3 \log \sigma^{-1}} \textrm{ and }
d_{k+1} - d_k \ll \alpha_k^{-2}\log \sigma^{-1}.
\end{equation*}
\end{enumerate}
From these last expressions we conclude that
\begin{equation*}
\alpha_k^2 \geq \alpha^2(1+2^{-4})^k,
\end{equation*}
and hence, since $\alpha_k \leq 1$, the iteration terminates with
$k \ll \log \alpha^{-1}$. It follows that
\begin{equation*}
d_k \ll \sum_{k=0}^\infty{\alpha_k^{-2}\log \sigma^{-2}} \leq
\alpha^{-2}\log \sigma^{-1}\sum_{k=0}^{\infty}{(1+2^{-4})^{-k}}
\ll \alpha^{-2} \log \sigma^{-1},
\end{equation*}
and that
\begin{equation*}
\delta_k \gg \left(\frac{\alpha}{\log \sigma^{-1}}\right)^{C\log
\alpha^{-1}}
\end{equation*}
for some absolute constant $C>0$.

For the iteration to terminate we must have arrived in the first
case at some point, and hence either
\begin{equation}\label{conclusion}
\sigma^{-1} \gg \left(\frac{\alpha}{\log
\sigma^{-1}}\right)^{C\log \alpha^{-1}}N^{c\alpha^2 (\log
\sigma^{-1})^{-1}}
\end{equation}
for some absolute constants $C,c>0$ or there is an arithmetic
progression in $A_1+A_2$ of length $(4\sigma)^{-1}$. The result
follows on taking $\sigma^{-1}$ as large a possible whilst not
satisfying (\ref{conclusion}).
\end{proof}

\section{Proof of Theorem \ref{aaageneral}}

As before we begin with a brief overview of the argument in the
finite-field setting, which can be made to prove the following
result.

\begin{theorem}
Suppose that $G$ is a compact vector space over $\mathbb{F}_2$ and
$A \subset G$ has density $\alpha>0$. Then $A+A+A$ contains (up to
a null set) an affine subspace of codimension at most
$4\alpha^{-1}$.
\end{theorem}
The proof is driven by the following iteration lemma.
\begin{lemma}
\emph{(Model iteration lemma)} Suppose that $G$ is a compact vector
space over $\mathbb{F}_2$. Suppose that $A \subset G$ has density
$\alpha$. Then at least one of the following is true.
\begin{enumerate}
\item $A+A+A$ contains $G$ (up to a null set). \item (Density
increment) There is a subspace $V$ of $G$ such that
\begin{equation*}
\|\chi_{A} \ast \mu_V\|_\infty \geq \alpha(1+\alpha/2) \textrm{
and } \cod V \leq 1.
\end{equation*}
\end{enumerate}
\end{lemma}
The iteration lemma is not conceptually difficult; we sketch the
main ideas now. Write $f:=\chi_A \ast \chi_A \ast \chi_A$. If $f$
is never zero (except for a null set) then $A+A+A$ certainly
contains $G$ (up to a null set), otherwise $f(x)=0$ on a set of
positive measure so there is a value of $x$ for which
\begin{equation*}
\sum_{\gamma \in \wh{G}}{\wh{\chi_A}(\gamma)^3\gamma(x)}=f(x)=0
\end{equation*}
by the inversion formula. Plancherel's theorem and the triangle
inequality in the usual fashion give
\begin{equation*}
\sum_{\gamma \neq 0_{\wh{G}}}{|\wh{\chi_A}(\gamma)|^3} \geq
\alpha^3,
\end{equation*}
from which it follows that there is a non-trivial characters
$\gamma$ at which $|\wh{f}(\gamma)|$ is large. Such a bound
provides an $\ell^\infty$ density increment for $A$.

Having proved this lemma the iteration is simple. Either $A+A+A$
contains a large affine subspace or we can increment the density
of $\alpha$. The density can't be incremented indefinitely and so
eventually $A+A+A$ contains a large affine subspace.

To localize the iteration argument is not as easy as it appears.
For the case $m=3$ the argument is really just Bourgain's original
argument for Roth's theorem. A particularly good exposition of
this, due to Tao, can be found in \cite{TT1}. There is a second
exposition also due to Tao in \cite{TT2} which uses smoothed
measures in place of our $\beta$s. For the generalization to $m>3$
the arguments in \cite{TT1} appear insufficient; in particular the
third claim in the proof below requires a new approach, which it
turns out was also used in \cite{TT2}. However, this is all,
perhaps, best illustrated by simply following the proof.
\begin{lemma}\emph{(Iteration lemma)} Suppose that $G$ is a
compact Abelian group and $B(\Gamma,\delta)$ a regular Bohr set in
$G$. Suppose that $A_1,...,A_m \subset B(\Gamma,\delta)$. Write
$\alpha$ for the geometric mean of the densities of the sets
$A_1,...,A_m$ in $B(\Gamma,\delta)$. Then at least one of the
following is true.
\begin{enumerate}
\item There is a $\delta'$ regular for $\Gamma$ such that
\begin{equation*}
\delta' \gg \min_i\{\int{\chi_{A_i}d\beta}\}^2\delta/md
\end{equation*}
and $A_1+...+A_m$ contains a translate of $B(\Gamma,\delta')$ (up
to a null set). \item There is a set of characters $\Gamma'$ and a
$\delta''$ regular for $\Gamma'$ such that
\begin{equation*}
|\Gamma'| \leq |\Gamma| +1, \delta'' \gg
\min_i\{\int{\chi_{A_i}d\beta}\}^3\delta/md^2
\end{equation*}
and
\begin{equation*}
\left(\prod_{i=1}^m{\|\chi_{A_i} \ast
\beta_{\Gamma',\delta''}\|_\infty}\right)^{\frac{1}{m}} \geq
\alpha\left(1+\frac{\alpha^{\frac{1}{m-2}}}{2^8m}\right).
\end{equation*}
\end{enumerate}
\end{lemma}
As with Lemma \ref{itlem} the proof which follows is rather complex
with a number of sub-claims being necessary. To ease understanding
we relegate proofs of these technical results to the end. The proof
itself essentially splits up the situation into the various ways in
which we can arrive at a density increment and then the technical
lemmas deal provide the density increments in each case.
\begin{proof}
We may certainly assume that $\alpha>0$ since otherwise we are done
for trivial reasons. Let $\delta'$ be a constant, regular for
$\Gamma$, to be chosen later. We may certainly assume that $A_1$ and
$A_2$ have the largest densities on $B(\Gamma,\delta)$ and so it is
$A_3,...,A_m$ we choose to move to the narrower Bohr neighborhood
$B(\Gamma,\delta')$.
\begin{equation*}
\int{\chi_{A_i} \ast \beta'd\beta}=\int{\chi_{A_i} d(\beta \ast
\beta')}=\int{\chi_{A_i}d\beta} + O(d\delta'\delta^{-1}),
\end{equation*}
by Corollary \ref{contlemcor}. It follows by averaging that there
is some $x_i \in B(\Gamma,\delta)$ such that
\begin{equation}\label{error11}
\chi_{A_i} \ast \beta'(x_i) \geq \int{\chi_{A_i}d\beta} +
O(d\delta'\delta^{-1}).
\end{equation}
Without loss of generality we may assume that the $x_i$s are all
zero. Write
\begin{eqnarray*}
& & \alpha_k := \int{\chi_{A_k} d\beta} \textrm{ and } f_k:=
(\chi_{A_k} - \alpha_k) \chi_{B} \textrm{ for } 1 \leq k \leq 2,\\
& \textrm{ and } & \alpha_k := \int{\chi_{A_k} d\beta'} \textrm{
and } f_k:= (\chi_{A_k} - \alpha_k) \chi_{B'} \textrm{ for } 3
\leq k \leq m.
\end{eqnarray*}

Define
\begin{equation*}
S:=B' \setminus \supp \chi_{A_1}\chi_B \ast \chi_{A_2}d\beta \ast
\chi_{A_3}d\beta' \ast ... \ast \chi_{A_m}d\beta'
\end{equation*}
and write $\sigma$ for the density of $S$ in $B'$. Now
$A_1+...+A_m \supset \supp \chi_{A_1}\chi_B \ast \chi_{A_2}d\beta
\ast \chi_{A_3}d\beta' \ast ... \ast \chi_{A_m}d\beta'$ so it
follows that if $\sigma=0$ then we are in the first case of the
lemma. Hence we assume that $\sigma>0$. We investigate the natural
inner product
\begin{equation}\label{ip}
I:=\langle f_1 \ast f_2d\beta \ast f_3d\beta' \ast ... \ast
f_md\beta', \chi_S \rangle_{L^2(\beta')}.
\end{equation}
We can rewrite $f_1 \ast f_2d\beta \ast f_3d\beta' \ast ... \ast
f_md\beta'$ as
\begin{eqnarray}
\label{decomp} & & \chi_{A_1}\chi_B \ast \chi_{A_2}d\beta \ast
\chi_{A_3}d\beta' \ast ... \ast \chi_{A_m}d\beta'\\ \nonumber & &
- \alpha_1 \chi_{B}\ast \chi_{A_2}d\beta \ast \chi_{A_3}d\beta'
\ast ... \ast \chi_{A_m}d\beta'\\ \nonumber & & - \chi_{A_1}\chi_B
\ast \alpha_2\beta \ast \chi_{A_3}d\beta' \ast ... \ast
\chi_{A_m}d\beta'\\ \nonumber & & + \alpha_1\chi_B \ast
\alpha_2\beta \ast \chi_{A_3}d\beta' \ast ... \ast
\chi_{A_m}d\beta' \\ \nonumber & & - f_1 \ast f_2d\beta \ast
\alpha_3\beta' \ast \chi_{A_4}d\beta' \ast ... \ast \chi_{A_m}d\beta'\\
\nonumber & & - ... \\ \nonumber & & - f_1 \ast f_2d\beta \ast
f_3d\beta' \ast ... \ast f_{j-1}d\beta' \ast \alpha_j\beta' \ast
\chi_{A_{j+1}}d\beta' \ast ... \ast \chi_{A_m}d\beta'\\ \nonumber
& & - ...\\ \nonumber & &  - f_1 \ast f_2d\beta \ast f_3d\beta'
\ast ... \ast f_{m-1}d\beta'\ast \alpha_m\beta'.
\end{eqnarray}
There are three different types of term in this decomposition. The
first term is unique and we denote it by $Z$, the next three are
all of the same type and we denote them by $T_1,T_2$ and $T_3$.
Finally the remaining terms are all of the same type and for $3
\leq j \leq m$ we write
\begin{equation*}
S_j=-f_1 \ast f_2d\beta \ast f_3d\beta' \ast ... \ast
f_{j-1}d\beta' \ast \alpha_j\beta' \ast \chi_{A_{j+1}}d\beta' \ast
... \ast \chi_{A_m}d\beta'.
\end{equation*}
We have
\begin{eqnarray}
\label{decompip}I& =& \langle
Z,\chi_S\rangle_{L^2(\beta')}+\langle
T_1,\chi_S\rangle_{L^2(\beta')}+\langle
T_2,\chi_S\rangle_{L^2(\beta')}+\langle
T_3,\chi_S\rangle_{L^2(\beta')}\\ \nonumber & & +\langle
S_3,\chi_S\rangle_{L^2(\beta')}+...+\langle
S_m,\chi_S\rangle_{L^2(\beta')}.
\end{eqnarray}
Our objective now is to estimate the inner products on the right.

The first inner product is zero since $\chi_S$ is supported on the
relative complement of $Z$. The inner products $\langle
T_i,\chi_S\rangle_{L^2(\beta')}$ can all be estimated in the same
way using the following claim which is Lemma \ref{techlem3}.
\begin{claim}
Suppose that $f \in L^\infty(\beta)$. Then
\begin{eqnarray*}
\langle  f\ast \beta \ast \chi_{A_3}d\beta' \ast ... \ast
\chi_{A_m}d\beta', \chi_S \rangle_{L^2(\beta')} & = & \sigma
\alpha_3...\alpha_m \times\\ & & \left(\int{fd\beta} +
O(md\delta'\delta^{-1}\|f\|_\infty)\right).
\end{eqnarray*}
\end{claim}

Note that the $T_i$s can be rewritten as follows.
\begin{eqnarray*}
T_1 & = & -\alpha_1\chi_{A_2}\chi_B \ast \beta \ast
\chi_{A_3}d\beta'
\ast ... \ast \chi_{A_m}d\beta'\\
T_2 & = & -\alpha_2\chi_{A_1}\chi_B \ast \beta \ast
\chi_{A_3}d\beta'
\ast ... \ast \chi_{A_m}d\beta'\\
T_3 & = & \alpha_1\alpha_2\chi_{B} \ast \beta \ast
\chi_{A_3}d\beta' \ast ... \ast \chi_{A_m}d\beta'.
\end{eqnarray*}
Now apply the claim in each case to see that
\begin{eqnarray*}
\langle T_1, \chi_S \rangle_{L^2(\beta')} & = &
-\sigma\alpha_1...\alpha_m( 1 + O(md\delta'\delta^{-1}\alpha_2^{-1}))\\
\langle T_2, \chi_S \rangle_{L^2(\beta')} & = &
-\sigma\alpha_1...\alpha_m( 1 + O(md\delta'\delta^{-1}\alpha_1^{-1}))\\
\langle T_3, \chi_S \rangle_{L^2(\beta')} & = &
\sigma\alpha_1...\alpha_m( 1 + O(md\delta'\delta^{-1})).
\end{eqnarray*}
It follows that
\begin{equation*}
\langle T_1, \chi_S \rangle_{L^2(\beta')}+...+\langle T_3, \chi_S
\rangle_{L^2(\beta')}=-\sigma\alpha_1...\alpha_m( 1 +
O(md\delta'\delta^{-1}(\alpha_1^{-1}+\alpha_2^{-1}))),
\end{equation*}
and hence in (\ref{decompip}) we have
\begin{equation*}
I - \langle S_3, \chi_S\rangle_{L^2(\beta')}- ... - \langle S_m,
\chi_S\rangle_{L^2(\beta')}= -\sigma\alpha_1...\alpha_m( 1 +
O(md\delta'\delta^{-1}(\alpha_1^{-1}+\alpha_2^{-1}))).
\end{equation*}
It follows that there is a $\delta'_0 \gg \delta
\min\{\alpha_1,\alpha_2\}/md$ such that if $\delta' \leq
\delta'_0$ then the error term here is at most $1/2$. We assume
that $\delta' \leq \delta'_0$ so that by the triangle inequality
\begin{equation*}
|I| + |\langle S_3, \chi_S\rangle_{L^2(\beta')}|+ ... +|\langle
S_m, \chi_S\rangle_{L^2(\beta')}|\geq
\frac{\sigma\alpha_1...\alpha_m}{2}.
\end{equation*}
It follows by averaging that one of the following is true.
\begin{equation*}
|I| \geq \sigma\alpha_1...\alpha_m/4 \textrm{ or } |\langle
S_j,\chi_S \rangle_{L^2(\beta')}| \geq \frac{\sigma
\alpha_1...\alpha_m}{2^j} \textrm{ for some } 3 \leq j \leq m.
\end{equation*}
We have two claims which deal with the two cases: they are proved in
Lemmas \ref{techlem4} and \ref{techlem5}.
\begin{claim}
If $|I| \geq \sigma\alpha_1...\alpha_m/4$ then there is a $k$ with
$3 \leq k \leq m$, a set of characters $\Gamma'$ and a $\delta''
\leq \delta'$ regular for $\Gamma'$ such that
\begin{equation*}
|\Gamma'| \leq |\Gamma|+1, \delta'' \gg \frac{\alpha_k^2\delta'}{d}
\textrm{ and } \|\chi_{A_k} \ast \beta_{\Gamma',\delta''}\|_\infty
\geq \alpha_k(1+2^{-5}\alpha^{\frac{1}{m-2}}).
\end{equation*}
\end{claim}
\begin{claim}
If
\begin{equation*}
|\langle S_j,\chi_S \rangle_{L^2(\beta')}| \geq \frac{\sigma
\alpha_1...\alpha_m}{2^j}
\end{equation*}
then either
\begin{enumerate}
\item there is a $k$ with $3 \leq k \leq j-1$, a set of characters
$\Gamma'$ and a $\delta'' \leq \delta'$ regular for $\Gamma'$ such
that
\begin{equation*}
|\Gamma'| \leq |\Gamma|+1, \delta'' \gg \frac{\alpha_k\delta'}{d}
\textrm{ and } \|\chi_{A_k} \ast \beta_{\Gamma',\delta''}\|_\infty
\geq \alpha_k(1+2^{-4});
\end{equation*} \item or there is an $k$ with $1 \leq k \leq 2$ and a $\delta'' \leq \delta'$ regular for $\Gamma$ such that
\begin{equation*}
\delta'' \gg \min\{\frac{\alpha\delta'}{d},\frac{\alpha_k
\delta}{d}\} \textrm{ and } \|\chi_{A_k} \ast
\beta_{\Gamma,\delta''}\|_\infty \geq \alpha_k(1+2^{-7}).
\end{equation*}
\end{enumerate}
\end{claim}
Now it is just a matter of choosing $\delta'$ as large as possible
whilst ensuring that the errors are small. From the claims we are
guaranteed at least one of the following three outcomes.
\begin{enumerate}
\item There is a $k$ with $3 \leq k \leq m$, a set of characters
$\Gamma'$ and a $\delta'' \leq \delta'$ regular for $\Gamma'$ such
that
\begin{equation*}
|\Gamma'| \leq |\Gamma|+1, \delta'' \gg
\frac{\alpha_k^2\delta'}{d} \textrm{ and } \|\chi_{A_k} \ast
\beta_{\Gamma',\delta''}\|_\infty \geq
\alpha_k(1+2^{-5}\alpha^{\frac{1}{m-2}}).
\end{equation*} \item There is a $k$ with $3 \leq k \leq m-1$, a set of characters
$\Gamma'$ and a $\delta'' \leq \delta'$ regular for $\Gamma'$ such
that
\begin{equation*}
|\Gamma'| \leq |\Gamma|+1, \delta'' \gg \frac{\alpha_k\delta'}{d}
\textrm{ and } \|\chi_{A_k} \ast \beta_{\Gamma',\delta''}\|_\infty
\geq \alpha_k(1+2^{-4});
\end{equation*} \item There is a $k$ with $1 \leq k \leq 2$ and a $\delta'' \leq \delta'$ regular for $\Gamma$ such that
\begin{equation*}
\delta'' \gg \min\{\frac{\alpha\delta'}{d},\frac{\alpha_k
\delta}{d}\} \textrm{ and } \|\chi_{A_k} \ast
\beta_{\Gamma,\delta''}\|_\infty \geq \alpha_k(1+2^{-7}).
\end{equation*}
\end{enumerate}
This imples that there is a $k$ with $1 \leq k \leq m$, a set of
characters $\Gamma'$ and a $\delta''$ regular for $\Gamma'$ with
\begin{equation*}
|\Gamma'| \leq |\Gamma|+1\textrm{ and } \delta'' \gg
\frac{\min\{\alpha,\alpha_k^2\}\delta'}{d}
\end{equation*}
such that
\begin{equation*}
\|\chi_{A_k} \ast \beta_{\Gamma',\delta''}\|_\infty \geq
\alpha_k(1+2^{-7}\alpha^{\frac{1}{m-2}}).
\end{equation*}
Moreover
\begin{eqnarray*}
\|\chi_{A_i} \ast \beta_{\Gamma',\delta''}\|_\infty & \geq &
\|\chi_{A_i} \ast \beta_{\Gamma',\delta''}\|_{L^1(\beta)}\\
& =& \int{\chi_{A_i}d(\beta \ast \beta'')}\\ & = &
\int{\chi_{A_i}d\beta} + O(d\delta''\delta^{-1}) \textrm{ by
Corollary \ref{contlemcor}.}
\end{eqnarray*}
It follows that
\begin{equation*}
\prod_{i=1}^m{\|\chi_{A_k} \ast \beta_{\Gamma',\delta''}\|_\infty}
\geq (1+2^{-7}\alpha^{\frac{1}{m-2}})\prod_{i=1}^m{
\left(\int{\chi_{A_i}d\beta} + O(d\delta''\delta^{-1})\right)}
\end{equation*}
which in turn is at least
\begin{equation*}
(1+2^{-7}\alpha^{\frac{1}{m-2}})\alpha^m.\prod_{i=1}^m{\left(1+O\left(d\delta''\delta^{-1}\left(\int{\chi_{A_i}d\beta}\right)^{-1}\right)\right)}.
\end{equation*}
When we take $m$th roots the product can be estimated by
\begin{equation*}
\left(\prod_{i=1}^m{\left(1+O\left(d\delta''\delta^{-1}\int{\chi_{A_i}d\beta}^{-1}\right)\right)}\right)^{\frac{1}{m}}
= 1+O(d\delta''\delta^{-1}\min_i\{\int{\chi_{A_i}d\beta}\}^{-1}),
\end{equation*}
so there is a $\delta'''_0$
\begin{equation*}
\delta'''_0 \gg
\frac{\delta\min_i\{\int{\chi_{A_i}d\beta}\}\alpha^{\frac{1}{m-2}}}{d}
\end{equation*}
such that if $\delta'' \leq \delta'''_0$ then
\begin{equation*}
\prod_{i=1}^m{\|\chi_{A_k} \ast
\beta_{\Gamma',\delta''}\|_\infty}^{\frac{1}{m}} \geq
\left(1+\frac{\alpha^{\frac{1}{m-2}}}{2^8m}\right)\alpha.
\end{equation*}
Since $\delta'' \leq \delta'$ the conclusion of the lemma follows
on taking $\delta'$ (regular by Proposition \ref{ubreg}) as large
as possible subject to $\delta' \leq \delta'_0$ and $\delta' \leq
\delta'''_0$.
\end{proof}

We now address the technical lemma which we employed above.

\begin{lemma}\label{techlem3}
Suppose that $G$ is a compact Abelian group, $B(\Gamma,\delta)$ is a
regular Bohr set, $B(\Gamma,\delta')$ is a Bohr set, $f \in
L^\infty(\beta)$ and $A_1,...,A_m,S \subset B(\Gamma,\delta')$ are
sets with relative density $\alpha_1,...,\alpha_m$ and $\sigma$
respectively. Then
\begin{eqnarray*}
\langle  f\ast \beta \ast \chi_{A_3}d\beta' \ast ... \ast
\chi_{A_m}d\beta', \chi_S \rangle_{L^2(\beta')} & = & \sigma
\alpha_3...\alpha_m \times\\ & & \left(\int{fd\beta} +
O(md\delta'\delta^{-1}\|f\|_\infty)\right).
\end{eqnarray*}
\end{lemma}
\begin{proof}
We show that if $x \in B'$ then $f\ast \beta \ast \chi_{A_3}d\beta'
\ast ... \ast \chi_{A_m}d\beta'(x)$ is a constant plus a small
$L^\infty$-error. This leads directly to the desired conclusion.

By Corollary \ref{contlemcor} with $\mu=\chi_{A_3}d\beta' \ast ...
\ast \chi_{A_m}d\beta'$ we have
\begin{equation*}
\|\beta \ast \mu - \alpha_3...\alpha_m\beta\| =
O(md\delta'\delta^{-1}\alpha_3...\alpha_m),
\end{equation*}
since $\supp \mu \subset B(\Gamma,m\delta')$. It follows that
\begin{eqnarray*}
f\ast \beta \ast \chi_{A_3}d\beta_3 \ast ... \ast \chi_{A_m}d\beta_m
& = & \alpha_3...\alpha_mf \ast \beta +
O(\|f\|_\infty \|\beta \ast \mu - \alpha_3...\alpha_m\beta\|)\\
& = & \alpha_3...\alpha_mf \ast \beta +
O(md\delta'\delta^{-1}\|f\|_\infty\alpha_{3}...\alpha_m).
\end{eqnarray*}
If $x \in B'$ then by Corollary \ref{contlemcor2}
\begin{equation*}
f \ast \beta(x) = \int{fd\beta}+O(\|f\|_\infty d\delta'\delta^{-1}).
\end{equation*}
Combining these last two expressions we get
\begin{equation*}
f\ast \beta \ast \chi_{A_3}d\beta' \ast ... \ast
\chi_{A_m}d\beta'(x)=\int{fd\beta}\alpha_3...\alpha_m+O(\|f\|_\infty
md\delta'\delta^{-1}\alpha_{3}...\alpha_m).
\end{equation*}
The required estimate follows.
\end{proof}
\begin{lemma}\label{techlem4}
Suppose that $G$ is a compact Abelian group, $B(\Gamma,\delta)$ and
$B(\Gamma,\delta')$ are regular Bohr sets, $A_1,A_2 \subset
B(\Gamma,\delta)$ have relative density $\alpha_1$ and $\alpha_2$
respectively, and $A_3,...,A_m,S \subset B(\Gamma,\delta')$ have
relative density $\alpha_3,...,\alpha_m$ and $\sigma>0$
respectively. Write $f_i:=(\chi_{A_i}-\alpha_i)\chi_B$ for $1 \leq i
\leq 2$ and $f_i:=(\chi_{A_i}-\alpha_i)\chi_{B'}$ for $3 \leq i \leq
m$. If
\begin{equation*}
|\langle f_1 \ast f_2d\beta \ast f_3d\beta' \ast ... \ast
f_md\beta',\chi_S \rangle_{L^2(\beta')}| \geq
\sigma\alpha_1...\alpha_m/4
\end{equation*}
then there is a $k$ with $3 \leq k \leq m$, a set of characters
$\Gamma'$ and a $\delta'' \leq \delta'$ regular for $\Gamma'$ such
that
\begin{equation*}
|\Gamma'| \leq |\Gamma|+1, \delta'' \gg \frac{\alpha_k^2\delta'}{d}
\textrm{ and } \|\chi_{A_k} \ast \beta_{\Gamma',\delta''}\|_\infty
\geq \alpha_k(1+2^{-5}\alpha^{\frac{1}{m-2}}).
\end{equation*}
\end{lemma}
\begin{proof}
Write $I:=\langle f_1 \ast f_2d\beta \ast f_3d\beta' \ast ... \ast
f_md\beta',\chi_S \rangle_{L^2(\beta')}$. Plancherel's theorem tells
us that
\begin{equation*}
I = \sum_{\gamma
\in\wh{G}}{\wh{f_1}(\gamma)\wh{f_2d\beta}(\gamma)\wh{f_3d\beta'}(\gamma)
...\wh{f_md\beta'}(\gamma)\overline{\wh{\chi_Sd\beta'}(\gamma)}}.
\end{equation*}
Recalling that $|\wh{\chi_Sd\beta'}(\gamma)| \leq
\|\chi_S\|_{L^1(\beta')} = \sigma$, we may apply the triangle
inequality to get
\begin{equation*}
\sigma\sum_{\gamma
\in\wh{G}}{|\wh{f_1}(\gamma)\wh{f_2d\beta}(\gamma)\wh{f_3d\beta'}(\gamma)
...\wh{f_md\beta'}(\gamma)|} \geq
\frac{\sigma\alpha_1...\alpha_m}{4}.
\end{equation*}
Divide by $\sigma$ (which is possible since $\sigma>0$) to get
\begin{equation}\label{leverage here}
\sum_{\gamma
\in\wh{G}}{|\wh{f_1}(\gamma)\wh{f_2d\beta}(\gamma)\wh{f_3d\beta'}(\gamma)
...\wh{f_md\beta'}(\gamma)|} \geq \frac{\alpha_1...\alpha_m}{4}.
\end{equation}
By the Cauchy-Schwarz inequality and Plancherel's theorem we have
\begin{eqnarray*}
\sum_{\gamma \not \in
\mathcal{L}}{|\wh{f_1}(\gamma)||\wh{f_2d\beta}(\gamma)|} & \leq &
\left(\sum_{\gamma \in
\wh{G}}{|\wh{f_1}(\gamma)|^2}\right)^{\frac{1}{2}}\left(\sum_{\gamma
\in \wh{G}}{|\wh{f_2d\beta}(\gamma)|^2}\right)^{\frac{1}{2}}\\
& = & \|f_1\|_2\|f_2d\beta\|_2\\ & = &
\|f_1\|_{L^2(\beta)}\|f_2\|_{L^2(\beta)} \textrm{ since $\beta$ is
uniform on $B$, }\\& = &
\left(\alpha_1(1-\alpha_1)\alpha_2(1-\alpha_2)\right)^{\frac{1}{2}}
\leq \sqrt{\alpha_1\alpha_2}=\alpha,
\end{eqnarray*}
so applying the triangle inequality to (\ref{leverage here}) we
conclude that
\begin{equation*}
\sup_{\gamma \in \wh{G}}{|\wh{f_3d\beta'}(\gamma)|
...|\wh{f_md\beta'}(\gamma)|}\geq
\frac{\alpha\alpha_3...\alpha_m}{4}.
\end{equation*}
It follows that for some $k$ with $3 \leq k \leq m$ we have
\begin{equation*}
\sup_{\gamma \in \wh{G}}{|\wh{f_kd\beta'}(\gamma)|^{m-2}} \geq
\frac{\alpha_k^{m-2}\alpha}{4}.
\end{equation*}
Apply Lemma \ref{linfty density} to get the lemma.
\end{proof}

\begin{lemma}\label{techlem5}
Suppose that $G$ is a compact Abelian group, $B(\Gamma,\delta)$ and
$B(\Gamma,\delta')$ are regular Bohr sets, $A_1,A_2 \subset
B(\Gamma,\delta)$ have relative density $\alpha_1$ and $\alpha_2$
respectively, and $A_3,...,A_m,S \subset B(\Gamma,\delta')$ have
relative density $\alpha_3,...,\alpha_m$ and $\sigma>0$
respectively. Write $f_i:=(\chi_{A_i}-\alpha_i)\chi_B$ for $1 \leq i
\leq 2$, $f_i:=(\chi_{A_i}-\alpha_i)\chi_{B'}$ for $3 \leq i \leq
m$, and
\begin{equation*}
S_j:=-f_1 \ast f_2d\beta \ast f_3d\beta' \ast ... \ast
f_{j-1}d\beta' \ast \alpha_j\beta' \ast \chi_{A_{j+1}}d\beta' \ast
... \ast \chi_{A_m}d\beta'
\end{equation*}
for $3 \leq j \leq m$. If
\begin{equation*}
|\langle S_j,\chi_S \rangle_{L^2(\beta')}| \geq \frac{\sigma
\alpha_1...\alpha_m}{2^j}
\end{equation*}
then either
\begin{enumerate}
\item there is a $k$ with $3 \leq k \leq j-1$, a set of characters
$\Gamma'$ and a $\delta'' \leq \delta'$ regular for $\Gamma'$ such
that
\begin{equation*}
|\Gamma'| \leq |\Gamma|+1, \delta'' \gg \frac{\alpha_k\delta'}{d}
\textrm{ and } \|\chi_{A_k} \ast \beta_{\Gamma',\delta''}\|_\infty
\geq \alpha_k(1+2^{-4});
\end{equation*} \item or there is an $k$ with $1 \leq k \leq 2$ and a $\delta'' \leq \delta'$ regular for $\Gamma$ such that
\begin{equation*}
\delta'' \gg \min\{\frac{\alpha\delta'}{d},\frac{\alpha_k
\delta}{d}\} \textrm{ and } \|\chi_{A_k} \ast
\beta_{\Gamma,\delta''}\|_\infty \geq \alpha_k(1+2^{-7}).
\end{equation*}
\end{enumerate}
\end{lemma}
\begin{proof}
Plancherel's theorem gives
\begin{equation*}
\langle S_j,\chi_S \rangle_{L^2(\beta')} = \sum_{\gamma \in
\wh{G}}{\wh{S_j}(\gamma)\overline{\wh{\chi_Sd\beta'}(\gamma)}}.
\end{equation*}
Recalling that $|\wh{\chi_Sd\beta'}(\gamma)| \leq
\|\chi_S\|_{L^1(\beta')} = \sigma$, we may apply the triangle
inequality to get
\begin{equation*}
|\langle S_j,\chi_S \rangle_{L^2(\beta')}| \leq \sigma\sum_{\gamma
\in \wh{G}}{|\wh{S_j}(\gamma)|}.
\end{equation*}
If we now use the assumption on the magnitude of the inner product
and divide by $\sigma$ (which is possible since $\sigma>0$) we get
\begin{equation}\label{Fourier offset}
2^{-j}\alpha_1...\alpha_{m} \leq \sum_{\gamma \in
\wh{G}}{\left\{\begin{aligned}
|\wh{f_1}(\gamma)||\wh{f_2d\beta}(\gamma)||\wh{f_3d\beta'}(\gamma)|...
|\wh{f_{j-1}d\beta'}(\gamma)|\\ \times
|\alpha_j\wh{\beta'}(\gamma)||\wh{\chi_{A_{j+1}}d\beta'}(\gamma)|...
|\wh{\chi_{A_m}d\beta'}(\gamma)|
\end{aligned}\right\}}.
\end{equation}
First we note that $|\wh{\chi_{A_{k}}d\beta'}(\gamma)| \leq
\alpha_k$ for $j+1\leq k \leq m$. Second if there is some $k$ with
$3 \leq k \leq j-1$ such that $|\wh{f_kd\beta'}(\gamma)| \geq
\alpha_k/2$ then we may apply Lemma \ref{linfty density} to get the
density increment in the first case of the conclusion of the lemma.
Hence we assume that $|\wh{f_kd\beta'}(\gamma)| \leq \alpha_k/2$ for
all $k$ with $3 \leq k \leq j-1$. These two observations serve to
tell us that each summand in (\ref{Fourier offset}) is bounded above
by
\begin{equation*}
|\wh{f_1}(\gamma)||\wh{f_2d\beta}(\gamma)||\wh{\beta'}(\gamma)|
2^{-(j-3)}\alpha_3...\alpha_m.
\end{equation*}
Hence
\begin{equation}\label{sft}
2^{-3}\alpha_1\alpha_2 \leq \sum_{\gamma \in
\wh{G}}{|\wh{f_1}(\gamma)||\wh{f_2d\beta}(\gamma)||\wh{\beta'}(\gamma)|}.
\end{equation}
The characters at which $|\wh{\beta'}(\gamma)|$ is large make a
significant contribution to this sum, which we can see as follows.
Write
\begin{equation*}
\mathcal{L}:=\{\gamma \in \wh{G}:|\wh{\beta'}(\gamma)| \geq
2^{-4}\sqrt{\alpha_1\alpha_2}\}.
\end{equation*}
Then
\begin{equation}\label{smallchars}
\sum_{\gamma \not\in
\mathcal{L}}{|\wh{f_1}(\gamma)||\wh{f_2d\beta}(\gamma)||\wh{\beta'}(\gamma)|}
\leq 2^{-4}\sqrt{\alpha_1\alpha_2} \sum_{\gamma \in
\wh{G}}{|\wh{f_1}(\gamma)||\wh{f_2d\beta}(\gamma)|}.
\end{equation}
Now by the Cauchy-Schwarz inequality and Plancherel's theorem we
have
\begin{eqnarray*}
\sum_{\gamma \not \in
\mathcal{L}}{|\wh{f_1}(\gamma)||\wh{f_2d\beta}(\gamma)|} & \leq &
\left(\sum_{\gamma \in
\wh{G}}{|\wh{f_1}(\gamma)|^2}\right)^{\frac{1}{2}}\left(\sum_{\gamma
\in \wh{G}}{|\wh{f_2d\beta}(\gamma)|^2}\right)^{\frac{1}{2}}\\
& = & \|f_1\|_2\|f_2d\beta\|_2\\ & = &
\|f_1\|_{L^2(\beta)}\|f_2\|_{L^2(\beta)} \textrm{ since $\beta$ is
uniform on $B$, }\\& = &
\left(\alpha_1(1-\alpha_1)\alpha_2(1-\alpha_2)\right)^{\frac{1}{2}}
\leq \sqrt{\alpha_1\alpha_2}.
\end{eqnarray*}
We can use this in (\ref{smallchars}) to see that
\begin{equation*}
\sum_{\gamma \not\in
\mathcal{L}}{|\wh{f_1}(\gamma)||\wh{f_2d\beta}(\gamma)||\wh{\beta'}(\gamma)|}
\leq 2^{-4}\alpha_1\alpha_2
\end{equation*}
and hence by (\ref{sft}) that
\begin{equation*}
\sum_{\gamma \in
\mathcal{L}}{|\wh{f_1}(\gamma)||\wh{f_2d\beta}(\gamma)||\wh{\beta'}(\gamma)|}
\geq 2^{-4}\alpha_1\alpha_2.
\end{equation*}
Apply Lemma \ref{nestsupport} to get a $\delta''_0 \gg
\sqrt{\alpha_1\alpha_2}\delta'/d$ such that for all $\delta'' \leq
\delta''_0$
\begin{eqnarray*}
\mathcal{L} & \subset & \{ \gamma \in \wh{G}: |1 - \gamma(x)| \leq
1/2 \textrm{ for all } x \in B(\Gamma,\delta'')\}\\ & \subset &
\{\gamma \in\wh{G}: |\wh{\beta_{\Gamma,\delta''}}(\gamma)| \geq
1/2\}.
\end{eqnarray*}
Write $\beta''$ for $\beta_{\Gamma,\delta''}$. Now, if $\gamma \in
\mathcal{L}$ we have
\begin{equation*}
|\wh{\beta''}(\gamma)|^2  \geq |\wh{\beta'}(\gamma)|/4,
\end{equation*}
so
\begin{equation*}
\sum_{\gamma \in
\wh{G}}{|\wh{f_1}(\gamma)||\wh{f_2d\beta}(\gamma)||\wh{\beta''}(\gamma)|^2}
\geq 2^{-6}\alpha_1\alpha_2.
\end{equation*}
Now applying the Cauchy-Schwarz inequality and Plancherel's theorem
as before, we get that the sum on the left bounded above by
\begin{eqnarray*}
& & \left(\sum_{\gamma \in
\wh{G}}{|\wh{f_1}(\gamma)\wh{\beta''}(\gamma)|^2}\right)^{\frac{1}{2}}\left(\sum_{\gamma
\in \wh{G}}{|\wh{f_2d\beta}(\gamma)\wh{\beta''}(\gamma)|^2}\right)^{\frac{1}{2}}\\
& = & \|f_1\ast \beta''\|_2\|(f_2d\beta)\ast \beta''\|_{2}\\
& = & \langle f_1 \ast \beta'', (f_1d\beta) \ast
\beta''\rangle^{\frac{1}{2}}\langle f_2 \ast \beta'', (f_2d\beta)
\ast \beta''\rangle^{\frac{1}{2}} \textrm{ since $\beta$ is uniform
on $B(\Gamma,\delta)$}.
\end{eqnarray*}
Hence
\begin{equation*}
\langle f_1 \ast \beta'', (f_1d\beta) \ast
\beta''\rangle^{\frac{1}{2}}\langle f_2 \ast \beta'', (f_2d\beta)
\ast \beta''\rangle^{\frac{1}{2}} \geq 2^{-6}\alpha_1\alpha_2.
\end{equation*}
It follows that there is some $k$ with $1 \leq k \leq 2$ such that
\begin{equation*}
\langle f_k \ast \beta'', (f_kd\beta) \ast \beta''\rangle \geq
2^{-6}\alpha_k^2,
\end{equation*}
and applying Lemma \ref{l2 density} we get
\begin{equation*}
\|\chi_{A_k} \ast \beta''\|_\infty \geq \alpha_k(1+2^{-6}) +
O(d\delta''\delta^{-1}).
\end{equation*}
It follows from Proposition \ref{ubreg} that there is a choice of
$\delta''$ regular for $\Gamma$ such that
\begin{equation*}
\min\{\delta_0'',\delta'\} \geq \delta'' \gg \min\{\delta''_0,\delta
\alpha_k/d\} \textrm{ and } \|\chi_{A_k} \ast \beta''\|_\infty \geq
\alpha_k(1+2^{-7}).
\end{equation*}
This gives the second conclusion of the lemma once we note that
$\alpha_1\alpha_2 \geq \alpha^2$.
\end{proof}

It is a simple matter to, as before, iterate this lemma.

\begin{proof}
[Proof of Theorem \ref{aaageneral}.] We construct a sequence of
regular Bohr sets $B(\Gamma_k,\delta_k)$ iteratively. Write
\begin{equation}
\beta_k=\beta_{\Gamma_k,\delta_k}, d_k:=|\Gamma_k| \textrm{ and }
\alpha_k = \left(\prod_{i=1}^m{\|\chi_{A_i} \ast
\beta_k\|_\infty}\right)^{\frac{1}{m}}.
\end{equation}
We initialize the iteration with $\Gamma_0=\{0_{\wh{G}}\}$ and
$\delta_0\gg 1$ regular for $\Gamma_0$ by Proposition \ref{ubreg}.

Suppose that we are at stage $k$ of the iteration. Since
$B(\Gamma_k,\delta_k)$ has positive measure $\chi_{A_i} \ast
\beta_k$ is continuous and hence we make take $x_1,...,x_m$ such
that
\begin{equation*}
\chi_{A_i} \ast \beta_k(x_i) = \|\chi_{A_i} \ast \beta_k\|_\infty.
\end{equation*}
Now we apply the iteration lemma to the sets $(A_1-x_1)\cap
B(\Gamma_k,\delta_k),(A_2-x_2)\cap
B(\Gamma_k,\delta_k),...,(A_m-x_m)\cap B(\Gamma_k,\delta_k)$ and
the regular Bohr set $B(\Gamma_k,\delta_k)$.
\begin{enumerate}
\item Either $A_1+...+A_m$ contains (up to a null set) a translate
of a Bohr set $B(\Gamma_k,\delta_k')$ with $\delta_k' \gg
\alpha_k^{2m}\delta_k/md_k$. \item Or there is a regular Bohr set
$B(\Gamma_{k+1},\delta_{k+1})$ such that
\begin{equation*}
\alpha_{k+1} \geq
\alpha_k\left(1+\frac{\alpha_k^{\frac{1}{m-2}}}{2^8m}\right),
\delta_{k+1} \gg \frac{\alpha_k^{3m}\delta_k}{m d_k^2} \textrm{
and } d_{k+1} - d_k \leq 1.
\end{equation*}
\end{enumerate}
From these last expressions we conclude that after at most
$2^8m\alpha_k^{-\frac{1}{m-2}}$ iterations the density doubles and
so the iteration terminates and moreover
\begin{equation*}
d_k \leq \sum_{j=0}^{\log_2\alpha^{-1}}
{2^8m(2^j\alpha)^{-\frac{1}{m-2}}} \leq \sum_{j=0}^{\infty}
{2^8m(2^j\alpha)^{-\frac{1}{m-2}}} \ll m^2\alpha^{-\frac{1}{m-2}},
\end{equation*}
and hence
\begin{equation*}
\delta_k \gg \left(c\alpha\right)^{Cm^3 \alpha^{-\frac{1}{m-2}}}
\end{equation*}
for some absolute constants $C,c>0$.

For the iteration to terminate we must have arrived at the first
case at some point and the conclusion follows.
\end{proof}

\section*{Acknowledgments}
I should like to thank Tim Gowers and Ben Green for supervision and
reading the drafts of this paper and an anonymous referee for some
very useful suggestions for clarification.

\bibliographystyle{alpha}

\bibliography{master}

\end{document}